\documentclass{amsart}
\usepackage[latin1]{inputenc}
\usepackage[english]{babel}
\usepackage{amsmath}
\usepackage{amsfonts}
\usepackage{amssymb}
\usepackage{latexsym}
\usepackage{enumerate}
\usepackage{pstricks-add}
\usepackage{graphicx}
\usepackage{ctable}
\usepackage{varioref}
\usepackage[vario]{fancyref}
\usepackage{xcolor,colortbl}
\usepackage{subfig}
\usepackage{physics}
\usepackage[section]{placeins}
\usepackage{mathtools}
\usepackage{qtree}
\usepackage{url}
\usepackage{hyperref}
\makeindex
\usepackage{enumerate}
\usepackage{multicol}
\newcommand{\transp}[1]{{#1}^\top}
\newcommand\numberthis{\addtocounter{equation}{1}\tag{\theequation}}
\newcommand{\cda}[1]{\mathbb{A}_{#1}}
\newcommand{\Gr}{W}
\newcommand{\x}{$\boxtimes$}
\newcommand{\twist}{\omega}

\numberwithin{equation}{section}
\numberwithin{table}{section}
\author{John W. Bales}
\address{Department of Mathematics(Retired)\\Tuskegee University\\Tuskegee, AL  36088}
\email{john.w.bales@gmail.com}
\curraddr{PO Box 210, Waverly, AL 36879}
\date{}
\title{The Cayley-Dickson doubling products}

\begin{document}

\begin{abstract}
The purpose of this paper is to identify all of the basic Cayley-Dickson doubling products. A Cayley-Dickson algebra $\cda{N+1}$ of dimension $2^{N+1}$ consists of all ordered pairs of elements of a Cayley-Dickson algebra $\cda{N}$ of dimension $2^N$ where the product $(a,b)(c,d)$ of elements of $\cda{N+1}$ is defined in terms of a pair of second degree binomials $\left(f(a,b,c,d),g(a,b,c,d)\right)$ satisfying certain properties. The polynomial pair$(f,g)$ is called a `doubling product.' While $\cda{0}$ may denote any ring, here it is taken to be the set $\mathbb{R}$ of real numbers. The binomials $f$ and $g$ should be devised such that $\cda{1}=\mathbb{C}$ the complex numbers, $\cda{2}=\mathbb{H}$ the quaternions, and $\cda{ 3}=\mathbb{O}$ the octonions. Historically, various researchers have used different yet equivalent doubling products.\\[4pt]
\noindent\textbf{Notice:} The published version of this paper is \copyright Springer and is available on their website at the following link:\\ \href{url}{https://link.springer.com/article/10.1007/s00006-015-0638-6}
\end{abstract}

\maketitle

\section{Introduction}

Although there are several possible Cayley-Dickson doubling products, past researchers have restricted themselves to only two of them. Furthermore, the development of the basis vectors has proceeded historically in a fashion which can obscure certain periodicities of the structure constants associated with the products of the basis vectors. The purpose of this paper is to catalog all possible variants of the Cayley-Dickson doubling products and to recommend a different way to number the basis vectors. The alternate numbering method has been used in the past, for example in Shafer's 1954 paper \cite{S1954}, but has fallen out of favor. By identifying the ordered pair of two sequences as the `shuffling' of the two sequences the author demonstrates how this leads naturally to the numbering used by Shafer and others and how it illustrates the periodicity of the structure constants and corresponding `twist' on the group product underlying the product of the basis vectors. As a result of cataloging different doubling products the author has identified one in particular $(a,b)(c,d)=(ac-b^*d ,da^*+bc )$ having an interesting twist map (although it does not satisfy the Moufang identities).
\section{Defining properties}\label{S:definingProperties}
We shall look for all doubling products satisfying the following necessary (but not sufficient) properties of a Cayley-Dickson algebra $\cda{N}$:
\begin{enumerate}
 \item For each $N$ and each $x\in \cda{N}$,
 \begin{align}\label{eq:unit}
   1\cdot x&=x\cdot1=x
 \end{align}
 
 \item There shall be a norm and an involution $^*$ such that for each $N$ and each $x,y\in \cda{N}$,
  
    \begin{align}
      \label{eq:conjugateSum}x+x^*&\in \cda{0}=\mathbb{R}\\
      \label{eq:conjugateProduct}xx^*&=x^*x=\Vert x\Vert^2\\
      \label{eq:productConjugate}(xy)^*&=y^*x^*
  \end{align}
 \item There shall be an infinite sequence $e_0,e_1,e_2,\cdots$ of unit vectors such that
  \begin{enumerate}
   \item For each $N\ge0$ if $p<2^N$ then $e_p$ is a basis vector for $\cda{N}$
   \item The set $\{\pm e_p \,\vert\, 0\le p < 2^N \}$ is a group under the product on $\cda{N}$.
   \item $e_0=1$ and if $p>0$ then 
   \begin{align}
   e_p^2 &= -e_0\label{eq:vectorSquare}
   \end{align}
   \item If $0\ne p\ne q\ne 0$ then the \emph{anti-symmetric} property holds:
   \begin{equation}\label{eq:symmetric}
   e_pe_q+e_qe_p =0
   \end{equation}
   \item  If $0\ne p\ne q\ne 0$ then there is an $r\ne0$ such that either $e_pe_q=e_r$ or $e_pe_q=-e_r$ and the \emph{quaternion property} holds:
   
   \begin{equation}\label{eq: quaternic}
    e_pe_q = e_r \text{ implies } e_qe_r = e_p
   \end{equation}  
  \end{enumerate}
\end{enumerate}

We will identify eight doubling products which satisfy these properties, four of which produce the octonions at the third doubling and four which do not.

\section{The basis vectors}
The entire approach taken in this paper is highly dependent upon the manner in which the basis vectors are chosen. They are deliberately chosen to be the basis vectors of the Hilbert space $\ell^2$ of square-summable sequences. That is, 
\begin{align*}
 e_0 &= 1,0,0,0,\cdots\\
 e_1 &= 0,1,0,0,\cdots\\
 e_2 &= 0,0,1,0,\cdots\\
 \vdots &
\end{align*}
In order to relate these basis vectors to the idea of ordered pairs, we identify a real number $r$ with the sequence
\begin{equation}
r = r,0,0,0,\cdots
\end{equation}
and the ordered pair of two sequences 
\begin{align*}
x &= x_0,x_1,x_2,\cdots\\
y &= y_0,y_1,y_2,\cdots
\end{align*}
with the `shuffle' of those two sequences
\begin{align*}
(x,y) &= x_0,y_0,x_1,y_1,x_2,y_2,\cdots
\end{align*}

This `shuffle basis' immediately leads to two results: an inductive development of the basis vectors and the proper definition of the involution.

These basis vectors defined inductively as ordered pairs are
\begin{align}\label{D: basisVectors}
 e_0      &= 1\\
 e_{2p}   &= (e_p,0)\\
 e_{2p+1} &= (0,e_p)
\end{align}

Furthermore, since $x+x^*\in\mathbb{R}$ is a requirement it is sufficient that

\begin{equation}
 x^* = x_0,-x_1,-x_2,\cdots
\end{equation}
with the result that
\begin{align*}
 (x,y)^* &= x_0,-y_0,-x_1,-y_1,\cdots
\end{align*}
which, expressed in terms of ordered pairs is
\begin{equation}
(x,y)^* = (x^*,-y)
\end{equation}
which is the traditional involution for Cayley-Dickson algebras. This also guarantees by induction that for each $\cda{N}$

\begin{equation}
 x^{**}=x
\end{equation}

\section{Devising an adequate product}

It is required that $1\cdot x=x\cdot 1=x$.

Let $(a,b)(c,d)=(f(a,b,c,d),g(a,b,c,d))$.

 For it to be the case that $(a,b)(1,0)=(a,b)$ and $(1,0)(c,d)=(c,d)$ it is sufficient that $f$ contain one of the terms $ac$ or $ca$ plus or minus some product of $b$ or $b^*$ with $d$ or $d^*$.

 Furthermore, it is sufficient that $g$ contain some product of $a$ or $a^*$ with $d$ plus a product of $c$ or $c^*$ with $b$.\\
 
 It is required that $xx^*=x^*x=\Vert x\Vert$. Thus\\ $(a,b)(a^*,-b)=\left(\Vert a\Vert^2+\Vert b\Vert^2,0\right)$ and 
      $(c^*,-d)(c,d)=\left(\Vert c\Vert^2+\Vert d\Vert^2,0\right)$.

 Thus $f$ must contain a product of $a$ and $c$ or $a^*$ and $c^*$ minus a product of either $d$ and $b^*$ or a product of $b$ and $d^*$. Combined with the previous result it suffices for $f$ to equal one of the following binomials:
 
 \begin{align*}
  f_0(a,b,c,d) &= ca-b^*d\\
  f_1(a,b,c,d) &= ca-db^*\\
  f_2(a,b,c,d) &= ac-b^*d\\
  f_ 3(a,b,c,d) &= ac-db^*\\
  f_4(a,b,c,d) &= ca-bd^*\\
  f_5(a,b,c,d) &= ca-d^*b\\
  f_6(a,b,c,d) &= ac-bd^*\\
  f_7(a,b,c,d) &= ac-d^*b
 \end{align*}
 
 We already know that $g$ must contain one of the following terms: $ad$, $da$, $a^*d$ or $da^*$.
 In order for $g$ to be zero when an ordered pair is multiplied by its conjugate it must be the case that, for each of these four options $g$ must contain, respectively, the  terms $c^*b$, $bc^*$, $cb$ and $bc$.

 Therefore $g$ must be one of the following four binomials:
 
 \begin{align*}
 g_0(a,b,c,d) &= da^*+bc\\
 g_1(a,b,c,d) &= a^*d+cb\\
 g_2(a,b,c,d) &= ad+c^*b\\
 g_ 3(a,b,c,d) &= da+bc^*\\
 \end{align*}
 
 This results in a combination of  32 potential Cayley-Dickson doubling products. And indeed all  32 can be shown to satisfy properties \eqref{eq:unit} through \eqref{eq:symmetric} above. It is the quaternion property \eqref{eq: quaternic} which reduces the number from thirty-two to eight. Four of these eight will produce, at the third doubling, algebras isomorphic to the octonions and four `pseudo' Cayley-Dickson products will not.

\section{Twisted Group Algebra}

The set $\Gr=\{0,1,2, 3,\cdots\}$ of \emph{whole numbers} is a group under the bit-wise `exclusive or' (XOR) operation on their binary representations, with group identity 0. For example, $27 \oplus 14 = 11011_b \oplus 01110_b =10101_b=21$. Rather than representing the group operation as addition, we will use juxtaposition. Thus for $p,q\in\Gr$ we will denote  $p\oplus q$ as simply $pq$. This group operation on $\Gr^2$ is pertinent since for $p,q\in\Gr$ and each of the  32 doubling products $(f_i,g_j)$ there is an $\twist_{ij}(p,q)\in\{-1,1\}$ such that
\begin{equation}\label{eq: baseProduct}
 e_pe_q=\twist_{ij}(p,q)e_{pq}
\end{equation}

The map $\twist_{ij}:\Gr\times\Gr\rightarrow\{-1,1\}$ is called a \emph{twist} on the group $\Gr$ making each of the algebras resulting from the  32 doubling products a \emph{twisted group algebra}.\cite{EL1969} \cite{BS1970} \cite{RW1971}.

\section{Interior points of $\Gr^2$ and their successors}

Property \eqref{eq: quaternic} is a property of the \emph{interior points} of $\Gr^2$.  The interior points of $\Gr^2$ are the points $(p,q)\in\Gr^2$ such that $0\ne p\ne q\ne0$. Every point of $\Gr^2$ has four \emph{successor} points $(2p,2q),(2p,2q+1),(2p+1,2q)(2p+1,2q+1)$. Interior points of $\Gr^2$ which are not a successor point of an interior point are \emph{initial} interior points.
To establish a general property of the interior points of $\Gr^2$ by induction one must first establish that the property holds for the initial interior points (the basis step). The inductive step consists of showing that if the property holds for an interior point $(p,q)$ then the property holds for the four successors of $(p,q)$.

 We will find that for the \emph{initial} interior points property \eqref{eq: quaternic} holds for only \emph{half} of the  32 possible doubling products. Then we will find that only 8 of those 16 doubling products will satisfy property \eqref{eq: quaternic} for the remaining interior points of $\Gr^2.$

Before beginning with the induction, let us re-express property \eqref{eq: quaternic} with property \eqref{eq: baseProduct} in mind.

\begin{equation}\label{eq: quaternic2}
 \text{If } 0\ne p\ne q\ne 0 \text{ and if } e_pe_q=e_{pq} \text{ then } e_qe_{pq}=e_p\text{ and }e_{pq}e_p=e_q.
\end{equation}

We will adopt the notation $(p,q,r)$ to mean that $e_pe_q=e_r$.  Since $e_pe_q+e_qe_p=0$ for $0\ne p\ne q$, one but not both of $(p,q,pq)$ or $(q,p,pq)$ must be the case. 

\section{The basis step of the inductive proof of property \eqref{eq: quaternic2}}

First we will determine which of the  32 doubling products satisfy property \eqref{eq: quaternic2} for basis vector products $e_pe_q$ where $(p,q)$ is an initial interior point of $\Gr^2$. The initial interior points are precisely those interior points $(p,q)$ for which either $p=1$ or $q=1$ or for which $p$ and $q$ differ by exactly 1. Put another way, the initial interior points of $\Gr^2$ consist of the ordered pairs of the forms $(2s,1),(1,2s),(2s+1,1),(1,2s+1),(2s,2s+1),(2s+1,2s)$ where $s>0$.

For each interior point of $\Gr^2$ we must have either one or the other but not both of the following quaternion properties:

\begin{equation}\label{eq: option1}
 (2s,1,2s+1)\text{ and }(1,2s+1,2s)\text{ and }(2s+1,2s,1)\tag{$Q$}
\end{equation}

\begin{equation}\label{eq: option2}
 (1,2s,2s+1)\text{ and }(2s,2s+1,1)\text{ and }(2s+1,1,2s)\tag{$\transp{Q}$}
\end{equation}

Consider the following

\begin{equation}\label{eq: initialOne}
e_{2s}e_1=(e_s,0)(0,1)=\begin{cases}(0,\phantom{-}e_s)=\phantom{-}e_{2s+1}  &\text{for }g_2,g_ 3\text{ with all }f_i\\
                                            (0,-e_s)=-e_{2s+1} &\text{for }g_0,g_1\text{ with all }f_i
                               \end{cases}
\end{equation}

\begin{equation}\label{eq: initialTwo}
e_{1}e_{2s+1}=(0,1)(0,e_s)=\begin{cases}(\phantom{-}e_s,0)=\phantom{-}e_{2s} &\text{for }f_4,f_5,f_6,f_7\text{ with all }g_j\\
                                        (-e_s,0)=-e_{2s} &\text{for }f_0,f_1,f_2,f_ 3\text{ with all }g_j
                               \end{cases}
\end{equation}

\begin{equation}\label{eq: initialThree}
e_{2s+1}e_{2s}=(0,e_s)(e_s,0)=\begin{cases}(0,\phantom{-}1)=\phantom{1}e_{1}  &\text{for }g_2,g_ 3\text{ with all } f_i\\
                                           (0,-1)=-e_{1} &\text{for }g_0,g_1\text{ with all } f_i
                               \end{cases}
\end{equation}

To make sense of these results we summarize in table \ref{tab: firstElim} which of the  32 doubling products satisfy a condition from either option \ref{eq: option1} or \ref{eq: option2}. Any doubling product (shown in gray) which satisfies a condition from both options fails property \eqref{eq: quaternic2} and therefore is not a valid Cayley-Dickson doubling product.

\tiny
\begin{table}
\caption{First Elimination of Doubling Products}\label{tab: firstElim}
\begin{tabular}{ccccccccc}
 & $f_0$ & $f_1$ & $f_2$ & $f_ 3$ & $f_4$ & $f_5$ & $f_6$ & $f_7$\\\toprule
$g_0$ & \ref{eq: option2}  & \ref{eq: option2}  & \ref{eq: option2}  & \ref{eq: option2}  & \cellcolor{lightgray}\ref{eq: option2}  & \cellcolor{lightgray}\ref{eq: option2}  & \cellcolor{lightgray}\ref{eq: option2}  & \cellcolor{lightgray}\ref{eq: option2} \\
      & \ref{eq: option2}  & \ref{eq: option2}  & \ref{eq: option2}  & \ref{eq: option2}  & \cellcolor{lightgray}\ref{eq: option1}  & \cellcolor{lightgray}\ref{eq: option1}  & \cellcolor{lightgray}\ref{eq: option1}  & \cellcolor{lightgray}\ref{eq: option1} \\\midrule
$g_1$ & \ref{eq: option2}  & \ref{eq: option2}  & \ref{eq: option2}  & \ref{eq: option2}  & \cellcolor{lightgray}\ref{eq: option2}  & \cellcolor{lightgray}\ref{eq: option2}  & \cellcolor{lightgray}\ref{eq: option2}  & \cellcolor{lightgray}\ref{eq: option2} \\
      & \ref{eq: option2}  & \ref{eq: option2}  & \ref{eq: option2}  & \ref{eq: option2}  & \cellcolor{lightgray}\ref{eq: option1}  & \cellcolor{lightgray}\ref{eq: option1}  & \cellcolor{lightgray}\ref{eq: option1}  & \cellcolor{lightgray}\ref{eq: option1} \\\midrule
$g_2$ & \cellcolor{lightgray}\ref{eq: option1}  & \cellcolor{lightgray}\ref{eq: option1}  & \cellcolor{lightgray}\ref{eq: option1}  & \cellcolor{lightgray}\ref{eq: option1}  & \ref{eq: option1}  & \ref{eq: option1}  & \ref{eq: option1}  & \ref{eq: option1} \\
      & \cellcolor{lightgray}\ref{eq: option2}  & \cellcolor{lightgray}\ref{eq: option2}  & \cellcolor{lightgray}\ref{eq: option2}  & \cellcolor{lightgray}\ref{eq: option2}  & \ref{eq: option1}  & \ref{eq: option1}  & \ref{eq: option1}  & \ref{eq: option1} \\\midrule
$g_ 3$ & \cellcolor{lightgray}\ref{eq: option1}  & \cellcolor{lightgray}\ref{eq: option1}  & \cellcolor{lightgray}\ref{eq: option1}  & \cellcolor{lightgray}\ref{eq: option1}  & \ref{eq: option1}  & \ref{eq: option1}  & \ref{eq: option1}  & \ref{eq: option1} \\
      & \cellcolor{lightgray}\ref{eq: option2}  & \cellcolor{lightgray}\ref{eq: option2}  & \cellcolor{lightgray}\ref{eq: option2}  & \cellcolor{lightgray}\ref{eq: option2}  & \ref{eq: option1}  & \ref{eq: option1}  & \ref{eq: option1}  & \ref{eq: option1} \\\bottomrule
\end{tabular}
\end{table}
\normalsize

\section{The inductive step of the proof of property \eqref{eq: quaternic2}}

Suppose $(p,q)$ is an interior point of $\Gr^2$ and that $(p,q,pq)$ implies that $(q,pq,p)$ and $(pq,p,q)$. For which of the sixteen remaining doubling products does property \eqref{eq: quaternic2} follow for all four successors of $(p,q)$?\\

For each of the following four pairs of triple conditions only conditions from one of each pair may follow from the conditions $(p,q,pq)$, $(q,pq,p)$, $(pq,p,q)$ and a valid doubling product.

\begin{enumerate}[I. ] 
\item Contradictory conditions $A$ and $\tilde{A}$
\begin{equation}\label{eq: A1}
 (2p,2q,2pq)\text{ and }(2q,2pq,2p)\text{ and }(2pq,2p,2q)\tag{$A$}
\end{equation}
or
\begin{equation}\label{eq: A2}
 (2q,2p,2pq)\text{ and }(2p,2pq,2q)\text{ and }(2pq,2q,2p)\tag{$\tilde{A}$}
\end{equation}

\item Contradictory conditions $B$ and $\tilde{B}$
\begin{equation}\label{eq: B1}
 (2p,2q+1,2pq+1)\text{ and }(2q+1,2pq+1,2p)\text{ and }(2pq+1,2p,2q+1)\tag{$B$}
\end{equation}
or
\begin{equation}\label{eq: B2}
 (2q+1,2p,2pq+1)\text{ and }(2p,2pq+1,2q+1)\text{ and }(2pq+1,2q+1,2p)\tag{$\tilde{B}$}
\end{equation}

\item Contradictory conditions $C$ and $\tilde{C}$
\begin{equation}\label{eq: C1}
 (2p+1,2q,2pq+1)\text{ and }(2q,2pq+1,2p+1)\text{ and }(2pq+1,2q,2q)\tag{$C$}
\end{equation}
or
\begin{equation}\label{eq: C2}
 (2q,2p+1,2pq+1)\text{ and }(2p+1,2pq+1,2q)\text{ and }(2pq+1,2q,2p+1)\tag{$\tilde{C}$}
\end{equation}

\item Contradictory conditions $D$ and $\tilde{D}$
\begin{equation}\label{eq: D1}
 (2p+1,2q+1,2pq)\text{ and }(2q+1,2pq,2p+1)\text{ and }(2pq,2p+1,2q+1)\tag{$D$}
\end{equation}
or
\begin{equation}\label{eq: D2}
 (2q+1,2p+1,2pq)\text{ and }(2p+1,2pq,2q+1)\text{ and }(2pq,2q+1,2p+1)\tag{$\tilde{D}$}
\end{equation}
\end{enumerate}

Suppose $(p,q)$ is an interior point of $\Gr^2$ and that $(p,q,pq)$, $(q,pq,p)$, $(pq,p,q)$.

\begin{enumerate}[1. ]
\item
\begin{align*}
e_{2p}e_{2q}&=(e_p,0)(e_q,0)\\
            &=\begin{cases}
                   (\phantom{-}e_pe_q,0)=(\phantom{-}e_{pq},0)=\phantom{-}e_{2pq}&\text{for }f_2,f_ 3,f_6,f_7\\
                             (-e_pe_q,0)=(          -e_{pq},0)=          -e_{2pq}&\text{for }f_0,f_1,f_4,f_5
              \end{cases} \numberthis \label{eq: inductive1A}
\end{align*}
So $(2p,2q,2pq)$ follows from conditions \ref{eq: A1} for $f_2$, $f_ 3$, $f_6$, $f_7$ and\\ $(2q,2p,2pq)$ follows from conditions \ref{eq: A2} for $f_0$, $f_1$, $f_4$, $f_5$.
\item
\begin{align*}
e_{2q}e_{2pq}&=(e_q,0)(e_{pq},0)\\
             &=\begin{cases}
                   (\phantom{-}e_qe_{pq},0)=(\phantom{-}e_p,0)=\phantom{-}e_{2p}&\text{for }f_2,f_ 3,f_6,f_7\\
                             (-e_qe_{pq},0)=          (-e_p,0)=          -e_{2p}&\text{for }f_0,f_1,f_4,f_5
               \end{cases} \numberthis \label{eq: inductive1B}
\end{align*}
So $(2q,2pq,2p)$ follows from conditions \ref{eq: A1} for $f_2$, $f_ 3$, $f_6$, $f_7$ and\\ $(2pq,2q,2p)$ follows from conditions \ref{eq: A2} for $f_0$, $f_1$, $f_4$, $f_5$.
\item
\begin{align*}
e_{2pq}e_{2p}&=(e_{pq},0)(e_p,0)\\
             &=\begin{cases}
                   (\phantom{-}e_{pq}e_p,0)=(\phantom{-}e_q,0)=\phantom{-}e_{2q}&\text{for }f_2,f_ 3,f_6,f_7\\
                             (-e_{pq}e_p,0)=          (-e_q,0)=          -e_{2q}&\text{for }f_0,f_1,f_4,f_5
               \end{cases} \numberthis  \label{eq: inductive1C}
\end{align*}
So $(2pq,2p,2q)$ follows from conditions \ref{eq: A1} for $f_2$, $f_ 3$, $f_6$, $f_7$ and\\ $(2p,2pq,2q)$ follows from conditions \ref{eq: A2} for $f_0$, $f_1$, $f_4$, $f_5$.
\item
\begin{align*}
e_{2p}e_{2q+1}&=(e_p,0)(0,e_q)\\
              &=\begin{cases}
                   (0,\phantom{-}e_pe_q)=(0,\phantom{-}e_{pq})=\phantom{-}e_{2pq+1}&\text{for }g_0,g_2\\
                   (0,          -e_pe_q)=(0,          -e_{pq})=          -e_{2pq+1}&\text{for }g_1,g_ 3
                \end{cases} \numberthis   \label{eq: inductive2A}
\end{align*}
So $(2p,2q+1,2pq+1)$ follows from conditions \ref{eq: B1} for $g_0,g_2$ and\\ $(2q+1,2p,2pq+1)$ follows from conditions \ref{eq: B2} for $g_1,g_ 3$.
\item
\begin{align*}
e_{2q+1}e_{2pq+1}&=(0,e_q)(0,e_{pq})\\
                 &=\begin{cases}
                   (\phantom{-}e_qe_{pq},0)=(\phantom{-}e_p,0)=\phantom{-}e_{2p}&\text{for }f_0,f_2,f_4,f_6\\
                   (          -e_qe_{pq},0)=(          -e_p,0)=           -e_{2p}&\text{for }f_1,f_ 3,f_5,f_7
                   \end{cases} \numberthis   \label{eq: inductive2B}
\end{align*}
So $(2q+1,2pq+1,2p)$ follows from conditions \ref{eq: B1} for $f_0,f_2,f_4,f_6$ and $(2pq+1,2q+1,2p)$ follows from conditions \ref{eq: B2} for $f_1,f_ 3,f_5,f_7$.
\item
\begin{align*}
e_{2pq+1}e_{2p}&=(0,e_{pq})(e_p,0)\\
               &=\begin{cases}
                   (0,\phantom{-}e_{pq}e_p)=(0,\phantom{-}e_q)=\phantom{-}e_{2q+1}&\text{for }g_0,g_2\\
                   (0,          -e_{pq}e_p)=(0,          -e_q)=          -e_{2q+1}&\text{for }g_1,g_ 3
                 \end{cases} \numberthis \label{eq: inductive2C}
\end{align*}
So $(2pq+1,2p,2q+1)$ follows from conditions \ref{eq: B1} for $g_0,g_2$ and\\ $(2p,2pq+1,2q+1)$ follows from conditions \ref{eq: B2} for $g_1,g_ 3$.
\item
\begin{align*}
e_{2p+1}e_{2q}&=(0,e_p)(e_q,0)\\
              &=\begin{cases}
                   (0,\phantom{-}e_pe_q)=(0,\phantom{-}e_{pq})=\phantom{-}e_{2pq}&\text{for }g_0,g_2\\
                   (0,          -e_pe_q)=(0,          -e_{pq})=          -e_{2pq}&\text{for }g_1,g_ 3
                \end{cases} \numberthis \label{eq: inductive 3A}
\end{align*}
So $(2pq+1,2p,2q+1)$ follows from conditions \ref{eq: C1} for $g_0,g_2$ and\\ $(2p,2pq+1,2q+1)$ follows from conditions \ref{eq: C2} for $g_1,g_ 3$.
\item
\begin{align*}
e_{2p}e_{2q+1}&=(e_p,0)(0,e_q)\\
              &=\begin{cases}
                   (0,\phantom{-}e_pe_q)=(0,\phantom{-}e_{pq})=\phantom{-}e_{2pq+1}&\text{for }g_0,g_2\\
                   (0,          -e_pe_q)=(0,          -e_{pq})=          -e_{2pq+1}&\text{for }g_1,g_ 3
                \end{cases} \numberthis \label{eq: inductive 3B}
\end{align*}
So $(2p,2q+1,2pq+1)$ follows from conditions \ref{eq: C1} for $g_0,g_2$ and\\ $(2q+1,2p,2pq+1)$ follows from conditions \ref{eq: C2} for $g_1,g_ 3$.
\item
\begin{align*}
e_{2q+1}e_{2pq+1}&=(0,e_q)(0,e_{pq})\\
                 &=\begin{cases}
                   (\phantom{-}e_qe_{pq},0)=(\phantom{-}e_p,0)=\phantom{-}e_{2p}&\text{for }f_0,f_2,f_4,f_6\\
                   (          -e_qe_{pq},0)=(          -e_p,0)=          -e_{2p}&\text{for }f_1,f_ 3,f_5,f_7
                   \end{cases} \numberthis \label{eq: inductive 3C}
\end{align*}
So $(2q+1,2pq+1,2p)$ follows from conditions \ref{eq: C1} for $f_0,f_2,f_4,f_6$ and $(2pq+1,2q+1,2p)$ follows from conditions \ref{eq: C2} for $f_1,f_ 3,f_5,f_7$.

\item
\begin{align*}
e_{2p+1}e_{2q+1}&=(0,e_p)(0,e_q)\\
                &=\begin{cases}
                   (\phantom{-}e_pe_q,0)=(\phantom{-}e_{pq},0)=\phantom{-}e_{2pq}&\text{for }f_0,f_2,f_4,f_6\\
                   (          -e_pe_q,0)=(          -e_{pq},0)=          -e_{2pq}&\text{for }f_1,f_ 3,f_5,f_7
                  \end{cases} \numberthis \label{eq: inductive4A}
\end{align*}
So $(2p+1,2q+1,2pq)$ follows from conditions \ref{eq: D1} for $f_0,f_2,f_4,f_6$ and $(2q+1,2p+1,2pq)$ follows from conditions \ref{eq: D2} for $f_1,f_ 3,f_5,f_7$.

\item
\begin{align*}
e_{2q+1}e_{2pq}&=(0,e_q)(e_{pq},0)\\
               &=\begin{cases}
                   (0,\phantom{-}e_qe_{pq})=(0,\phantom{-}e_p)=\phantom{-}e_{2p+1}&\text{for }g_0,g_2\\
                   (0,          -e_qe_{pq})=(0,          -e_p)=          -e_{2p+1}&\text{for }g_1,g_ 3
                 \end{cases} \numberthis  \label{eq: inductive4B}
\end{align*}
So $(2q+1,2pq,2p+1)$ follows from conditions \ref{eq: D1} for $g_0,g_2$ and\\ $(2pq,2q+1,2p+1)$ follows from conditions \ref{eq: D2} for $g_1,g_ 3$.

\item
\begin{align*}
e_{2pq}e_{2p+1}&=(pq,0)(0,e_p)\\
               &=\begin{cases}
                   (0,\phantom{-}e_{pq}e_p)=(0,\phantom{-}e_q)=\phantom{-}e_{2q+1}&\text{for }g_0,g_2\\
                   (0,          -e_{pq}e_p)=(0,          -e_q)=          -e_{2q+1}&\text{for }g_1,g_ 3
                 \end{cases} \numberthis \label{eq: inductive4C}
\end{align*}
So $(2pq,2p+1,2q+1)$ follows from conditions \ref{eq: D1} for $g_0,g_2$ and\\ $(2p+1,2pq,2q+1)$ follows from conditions \ref{eq: D2} for $g_1,g_ 3$.\\
\end{enumerate}

These results are compiled in table \ref{tab: secondElim} on page \pageref{tab: secondElim}. Eight more of the doubling products are eliminated for inconsistency with property \ref{eq: quaternic2}. Thus the inductive step of the proof of property \ref{eq: quaternic2} only succeeds for the eight remaining products.

\tiny
\begin{table}[ht]
\caption{Second Elimination of Doubling Products}\label{tab: secondElim}
\begin{tabular}{lccccccccc}
      &&$f_0$ &$f_1$ &$f_2$ &$f_3$ &$f_4$ &$f_5$ &$f_6$ &$f_7$\\\toprule
 $g_0$&\ref{eq: A1}& &\cellcolor{lightgray}&\x&\cellcolor{lightgray}\x&\cellcolor{lightgray}&\cellcolor{lightgray}&\cellcolor{lightgray}&\cellcolor{lightgray}\\
      &\ref{eq: A2}&\x&\cellcolor{lightgray}\x& &\cellcolor{lightgray}&\cellcolor{lightgray}&\cellcolor{lightgray}&\cellcolor{lightgray}&\cellcolor{lightgray}\\\cline{2-10}
      &\ref{eq: B1}&\x&\cellcolor{lightgray}\x&\x&\cellcolor{lightgray}\x&\cellcolor{lightgray}&\cellcolor{lightgray}&\cellcolor{lightgray}&\cellcolor{lightgray}\\
      &\ref{eq: B2}& &\cellcolor{lightgray}\x& &\cellcolor{lightgray}\x&\cellcolor{lightgray}&\cellcolor{lightgray}&\cellcolor{lightgray}&\cellcolor{lightgray}\\\cline{2-10}
      &\ref{eq: C1}&\x&\cellcolor{lightgray}\x&\x&\cellcolor{lightgray}\x&\cellcolor{lightgray}&\cellcolor{lightgray}&\cellcolor{lightgray}&\cellcolor{lightgray}\\
      &\ref{eq: C2}& &\cellcolor{lightgray}\x& &\cellcolor{lightgray}\x&\cellcolor{lightgray}&\cellcolor{lightgray}&\cellcolor{lightgray}&\cellcolor{lightgray}\\\cline{2-10}
      &\ref{eq: D1}&\x&\cellcolor{lightgray}\x&\x&\cellcolor{lightgray}\x&\cellcolor{lightgray}&\cellcolor{lightgray}&\cellcolor{lightgray}&\cellcolor{lightgray}\\
      &\ref{eq: D2}& &\cellcolor{lightgray}\x& &\cellcolor{lightgray}\x&\cellcolor{lightgray}&\cellcolor{lightgray}&\cellcolor{lightgray}&\cellcolor{lightgray}\\\midrule
 $g_1$&\ref{eq: A1}&\cellcolor{lightgray} & &\cellcolor{lightgray}\x&\x&\cellcolor{lightgray}&\cellcolor{lightgray}&\cellcolor{lightgray}&\cellcolor{lightgray}\\
      &\ref{eq: A2}&\cellcolor{lightgray}\x&\x&\cellcolor{lightgray} & &\cellcolor{lightgray}&\cellcolor{lightgray}&\cellcolor{lightgray}&\cellcolor{lightgray}\\\cline{2-10}
      &\ref{eq: B1}&\cellcolor{lightgray} & &\cellcolor{lightgray}\x& &\cellcolor{lightgray}&\cellcolor{lightgray}&\cellcolor{lightgray}&\cellcolor{lightgray}\\
      &\ref{eq: B2}&\cellcolor{lightgray}\x&\x&\cellcolor{lightgray}\x&\x&\cellcolor{lightgray}&\cellcolor{lightgray}&\cellcolor{lightgray}&\cellcolor{lightgray}\\\cline{2-10}
      &\ref{eq: C1}&\cellcolor{lightgray}\x& &\cellcolor{lightgray}\x& &\cellcolor{lightgray}&\cellcolor{lightgray}&\cellcolor{lightgray}&\cellcolor{lightgray}\\
      &\ref{eq: C2}&\cellcolor{lightgray}\x&\x&\cellcolor{lightgray}\x&\x&\cellcolor{lightgray}&\cellcolor{lightgray}&\cellcolor{lightgray}&\cellcolor{lightgray}\\\cline{2-10}
      &\ref{eq: D1}&\cellcolor{lightgray}\x& &\cellcolor{lightgray}\x& &\cellcolor{lightgray}&\cellcolor{lightgray}&\cellcolor{lightgray}&\cellcolor{lightgray}\\
      &\ref{eq: D2}&\cellcolor{lightgray}\x&\x&\cellcolor{lightgray}\x&\x&\cellcolor{lightgray}&\cellcolor{lightgray}&\cellcolor{lightgray}&\cellcolor{lightgray}\\\midrule
 $g_2$&\ref{eq: A1}&\cellcolor{lightgray}&\cellcolor{lightgray}&\cellcolor{lightgray}&\cellcolor{lightgray}& &\cellcolor{lightgray} &\x&\cellcolor{lightgray}\x\\
      &\ref{eq: A2}&\cellcolor{lightgray}&\cellcolor{lightgray}&\cellcolor{lightgray}&\cellcolor{lightgray}&\x&\cellcolor{lightgray}\x& &\cellcolor{lightgray} \\\cline{2-10}
      &\ref{eq: B1}&\cellcolor{lightgray}&\cellcolor{lightgray}&\cellcolor{lightgray}&\cellcolor{lightgray}&\x&\cellcolor{lightgray}\x&\x&\cellcolor{lightgray}\x\\
      &\ref{eq: B2}&\cellcolor{lightgray}&\cellcolor{lightgray}&\cellcolor{lightgray}&\cellcolor{lightgray}& &\cellcolor{lightgray}\x& &\cellcolor{lightgray}\x\\\cline{2-10}
      &\ref{eq: C1}&\cellcolor{lightgray}&\cellcolor{lightgray}&\cellcolor{lightgray}&\cellcolor{lightgray}&\x&\cellcolor{lightgray}\x&\x&\cellcolor{lightgray}\x\\
      &\ref{eq: C2}&\cellcolor{lightgray}&\cellcolor{lightgray}&\cellcolor{lightgray}&\cellcolor{lightgray}& &\cellcolor{lightgray}\x& &\cellcolor{lightgray}\x\\\cline{2-10}
      &\ref{eq: D1}&\cellcolor{lightgray}&\cellcolor{lightgray}&\cellcolor{lightgray}&\cellcolor{lightgray}&\x&\cellcolor{lightgray}\x&\x&\cellcolor{lightgray}\x\\
      &\ref{eq: D2}&\cellcolor{lightgray}&\cellcolor{lightgray}&\cellcolor{lightgray}&\cellcolor{lightgray}& &\cellcolor{lightgray}\x& &\cellcolor{lightgray}\x\\\midrule
 $g_3$&\ref{eq: A1}&\cellcolor{lightgray}&\cellcolor{lightgray}&\cellcolor{lightgray}&\cellcolor{lightgray}&\cellcolor{lightgray} & &\cellcolor{lightgray}\x&\x\\
      &\ref{eq: A2}&\cellcolor{lightgray}&\cellcolor{lightgray}&\cellcolor{lightgray}&\cellcolor{lightgray}&\cellcolor{lightgray}\x&\x&\cellcolor{lightgray} & \\\cline{2-10}
      &\ref{eq: B1}&\cellcolor{lightgray}&\cellcolor{lightgray}&\cellcolor{lightgray}&\cellcolor{lightgray}&\cellcolor{lightgray}\x& &\cellcolor{lightgray}\x& \\
      &\ref{eq: B2}&\cellcolor{lightgray}&\cellcolor{lightgray}&\cellcolor{lightgray}&\cellcolor{lightgray}&\cellcolor{lightgray}\x&\x&\cellcolor{lightgray}\x&\x\\\cline{2-10}
      &\ref{eq: C1}&\cellcolor{lightgray}&\cellcolor{lightgray}&\cellcolor{lightgray}&\cellcolor{lightgray}&\cellcolor{lightgray}\x& &\cellcolor{lightgray}\x& \\
      &\ref{eq: C2}&\cellcolor{lightgray}&\cellcolor{lightgray}&\cellcolor{lightgray}&\cellcolor{lightgray}&\cellcolor{lightgray}\x&\x&\cellcolor{lightgray}\x&\x\\\cline{2-10}
      &\ref{eq: D1}&\cellcolor{lightgray}&\cellcolor{lightgray}&\cellcolor{lightgray}&\cellcolor{lightgray}&\cellcolor{lightgray}\x& &\cellcolor{lightgray}\x& \\
      &\ref{eq: D2}&\cellcolor{lightgray}&\cellcolor{lightgray}&\cellcolor{lightgray}&\cellcolor{lightgray}&\cellcolor{lightgray}\x&\x&\cellcolor{lightgray}\x&\x\\\midrule
\bottomrule
\end{tabular}
\end{table}
\normalsize

\noindent      $P_0$: $(a,b)(c,d)=(ca-b^*d ,da^*+bc )$\\
      $P_1$: $(a,b)(c,d)=(ca-db^* ,a^*d+cb )$\\
      $P_2$: $(a,b)(c,d)=(ac-b^*d ,da^*+bc )$\\
      $P_3$: $(a,b)(c,d)=(ac-db^* ,a^*d+cb )$\\
      $\transp{P}_0$: $(a,b)(c,d)=(ca-bd^* ,ad+c^*b )$\\
      $\transp{P}_1$: $(a,b)(c,d)=(ca-d^*b ,da+bc^* )$\\
      $\transp{P}_2$: $(a,b)(c,d)=(ac-bd^* ,ad+c^*b )$\\
      $\transp{P}_3$: $(a,b)(c,d)=(ac-d^*b ,da+bc^* )$\\

The `transpose' symbol is used for the second set of four doubling products since, as will later become apparent,  the corresponding product matrices of unit vectors are transposes of each other.

The product $\transp{P}_3$ is the one most commonly used, but its transpose $P_3$ is a close second. The author has found no instance of the use of the other six doubling products.

All eight of these products result in the complex $i=e_1=0,1,0,0,\cdots$. For the first four products, the quaternion $j=e_2=0,0,1,0,0,\cdots$ and $k=e_3=0,0,0,1,0,\cdots$. For their transposes $k=e_2$ and $j=e_3$.

Only four of these products, $P_0,\transp{P}_0,P_3,\transp{P}_3$ produce, at the third doubling, algebras isomorphic to the octonions. The other four are interesting for a different reason as will be shown later by looking at product $P_2$.

\section{Recursive definition of structure constants for the doubling products}

For $P_{0},P_{1},P_{2}$ and $P_{3},$ 
 
\begin{equation}(1,2n,2n+1) \text{\ for all\ } n>0 \label{eq: inductStartPlain}
\end{equation}

whereas for  $\transp{P}_0,\transp{P}_1,\transp{P}_2$ and $\transp{P}_3$, 

\begin{equation}(1,2n+1,2n) \text{\ for all\ } n>0\label{eq: inductStartHat}
\end{equation}

The `transpositive' nature of these two properties induces a transpose relationship between pairs of the eight doubling products.

For all eight of the doubling products, the second quaternion property \ref{eq: quaternic2} holds.

For $0\ne p\ne q\ne 0$

\begin{equation}(p,q,r)\longrightarrow (q,r,p)\longrightarrow (r,p,q)\end{equation}

For $P_{0}$ and $\transp{P}_{0}$ if $0\ne p\ne q\ne 0$ then 
\begin{align}
      (p,q,r) &\longrightarrow (2r,2q,2p)\notag\\
                    &\longrightarrow (2p,2q+1,2r+1)\notag\\
                    &\longrightarrow (2p+1,2q,2r+1)\notag\\
                    &\longrightarrow (2p+1,2q+1,2r) \label{eq: induct0}
\end{align}

For $P_{1}$ and $\transp{P}_{1}$ if $0\ne p\ne q\ne 0$ then 
\begin{align}
      (p,q,r) &\longrightarrow (2r,2q,2p)\notag\\
                    &\longrightarrow (2r+1,2q+1,2p)\notag\\
                    &\longrightarrow (2r+1,2q,2p+1)\notag\\
                    &\longrightarrow (2r,2q+1,2p+1) \label{eq: induct1}
\end{align}

For $P_{2}$ and $\transp{P}_{2}$ if $0\ne p\ne q\ne 0$ then 
\begin{align}
      (p,q,r) &\longrightarrow (2p,2q,2r)\notag\\
                    &\longrightarrow (2p,2q+1,2r+1)\notag\\
                    &\longrightarrow (2p+1,2q,2r+1)\notag\\
                    &\longrightarrow (2p+1,2q+1,2r) \label{eq: induct2}
\end{align}

For $P_{3}$ and $\transp{P}_{3}$ if $0\ne p\ne q\ne 0$ then 
\begin{align}
      (p,q,r) &\longrightarrow (2p,2q,2r)\notag\\
                    &\longrightarrow (2r+1,2q+1,2p)\notag\\
                    &\longrightarrow (2r+1,2q,2p+1)\notag\\
                    &\longrightarrow (2r,2q+1,2p+1) \label{eq: induct3}
\end{align}

The most commonly used of these eight doubling products is $\transp{P}_3$ \cite{S1954} \cite{B1967} \cite{M1998} but $P_3$ has also been used \cite{B2001} \cite{B2009}.
 
\section{Application to octonion basis vectors}

Using induction rules \vref{eq: inductStartPlain} through \vref{eq: induct3} one can construct all the $(p,q,r)$ for any $\cda{N}$. First using either rule \ref{eq: inductStartPlain} or \vref{eq: inductStartHat} for all $n$ such that $2n+1<2^N$ construct all $(1,q,r)$. Then use whichever of rules \vref{eq: induct0} through \vref{eq: induct3} applies to compute the remainder.

For example, using ${P}_0$, let us compute all the quaternion triplets for the octonions $\mathbb{O}=\cda{3}$. Applying rule \vref{eq: inductStartHat} for $n=1,2,3$ gives the quaternion triplets $(1,2,3), (1,4,5)$ and $(1,6,7)$. The remaining four can be obtained from $(1,2,3)$ and rule \vref{eq: induct0}: $(2,6,4), (2,5,7), (6,3,5)$ and $(4,7,3)$.

The following are the triplets for $P_0$ through $P_3$.

\begin{enumerate}
 \item[$P_0$: ]$(1,2,3),(1,4,5),(1,6,7),(2,6,4),(7,2,5),(5,6,3),(3,4,7)$
 \item[$P_1$: ]$(1,2,3),(1,4,5),(1,6,7),(2,6,4),(5,2,7),(7,4,3),(3,6,5)$
 \item[$P_2$: ]$(1,2,3),(1,4,5),(1,6,7),(2,4,6),(7,2,5),(5,6,3),(3,4,7)$
 \item[$P_3$: ]$(1,2,3),(1,4,5),(1,6,7),(2,4,6),(5,2,7),(7,4,3),(3,6,5)$
\end{enumerate}

 Reverse these to obtain the triplets for $\transp{P}_0$ through $\transp{P}_3$.
 
\begin{enumerate}
 \item[$\transp{P}_0$: ]$(1,3,2),(1,5,4),(1,7,6),(2,4,6),(7,5,2),(5,3,6),(3,7,4)$
 \item[$\transp{P}_1$: ]$(1,3,2),(1,5,4),(1,7,6),(2,4,6),(5,7,2),(7,3,4),(3,5,6)$
 \item[$\transp{P}_2$: ]$(1,3,2),(1,5,4),(1,7,6),(2,6,4),(7,5,2),(5,3,6),(3,7,4)$
 \item[$\transp{P}_3$: ]$(1,3,2),(1,5,4),(1,7,6),(2,6,4),(5,7,2),(7,3,4),(3,5,6)$
\end{enumerate}
 
For $P_0,\,P_3$ there are permutations of integers 1 through 7 which take any triple $(p,q,r)$ through all the triples for that product as well as its transpose. For example:

\begin{enumerate}
\item[$P_0$:]$(1263457)\vert(1,2,3),(2,6,4),(6,3,5),(3,4,7),(4,5,1),(5,7,2),(7,1,6)$
\item[$P_3$:]$(1243675)\vert(1,2,3),(2,4,6),(4,3,7),(3,6,5),(6,7,1),(7,5,2),(5,1,4)$
\end{enumerate}

\section{The oriented Fano plane}

The products of the basis vectors for the octonions are customarily represented in the Fano Plane. The shuffle basis produces an especially nice representation of the eight doubling products. (See Figure \vref{fig: fano}) Fano planes are typically numbered in a haphazard fashion although a better oriented construction exists. To construct an oriented Fano plane draw an equilateral triangle inscribed with a circle and construct the three altitudes. Label the center $e_1$. Next label a midpoint of one of the sides $e_2$. Label the vertex opposite that midpoint $e_3$. Label one of the two remaining midpoints $e_4$ and the vertex opposite that midpoint $e_5$. Label the remaining midpoint $e_6$ and the remaining vertex $e_7$. This results in one of only two versions of an oriented Fano plane, each being a reflection of the other in one of the altitudes.

Each of the three sides of the triangle represents a triple $(p,q,r).$ Likewise each altitude from a vertex to the midpoint of the opposite side as well as the circle through the midpoints represents such a triple.

For each of the eight doubling products, the sense of the three sides is the same--either clockwise ($\leftarrow$) around the triangle, or counter-clockwise ($\rightarrow$). If clockwise, then the three sides of the triangle represent $(5,2,7),$ $(7,4,3)$ and $(3,6,5).$ If counter-clockwise then the three sides represent $(7,2,5),$ $(5,6,3)$ and $(3,4,7).$ 

The circle through the midpoints of the sides represents either $(2,4,6)$ in the clockwise sense ($\circlearrowright$) or $(2,6,4)$ in the counter-clockwise sense ($\circlearrowleft$).

The three altitudes may all be in an `up' sense from base to vertex ($\uparrow$) or may all be in a `down' sense from vertex to base ($\downarrow$). So the altitudes must either be $(1,2,3),$ $(1,4,5),$ and $(1,6,7)$ or $(1,3,2),$ $(1,5,4),$ and $(1,7,6)$. All altitudes must have the same sense. See Figure \vref{fig: fano} for a breakdown of all the modes of the eight doubling products.

This common orientation with regards to directions of the three sides as well as directions of the altitudes, together with the fact that for each $(p,q,r)$, $r$ is the `bit-wise exclusive or' of $p$ and $q$, motivates calling this an `oriented' Fano plane. The sides all have the same orientation and the altitudes all have the same orientation.

Since the sides may have two senses and the circle may have two senses and the altitudes may have two senses, all $2^3=8$ versions of the doubling products may be accommodated in the one diagram.

For any of the eight doubling product variations, knowing the `sense' of $(1,2,3),$ $(2,4,6)$ and $(2,5,7)$ is sufficient to recover the sense of all triples $(p,q,r)$ for each product using an oriented Fano plane.

\begin{figure}[ht]
 \psset{xunit=1.0cm,yunit=1.0cm,algebraic=true,dotstyle=o,dotsize=3pt 0,
linewidth=1.5pt,arrowsize=3pt 2,arrowinset=0.25}
\begin{pspicture*}(-5,-3.75)(5,6.25)
\pspolygon[fillstyle=solid,fillcolor=white](-4.33,-2.5)(0,5)(4.33,-2.5)
\psline(-4.33,-2.5)(2.165,1.25)
\psline(0,5)(0,-2.5)
\psline(4.33,-2.5)(-2.165,1.25)
\pscircle(0,0){2.5}
\psdots[dotstyle=*,dotsize=6pt,linecolor=black](-4.33,-2.5)
\rput[bl](-4.5,-3){\black{$\mathbf{e_7}$}}
\psdots[dotstyle=*,dotsize=6pt,linecolor=black](4.33,-2.5)
\rput[bl](4.2,-3){\black{$\mathbf{e_5}$}}
\psdots[dotstyle=*,dotsize=6pt,linecolor=black](0,5)
\rput[bl](-0.1,5.2){\black{$\mathbf{e_3}$}}
\psdots[dotstyle=*,dotsize=6pt,linecolor=black](2.165,1.25)
\rput[bl](2.2,1.35){\black{$\mathbf{e_6}$}}
\psdots[dotstyle=*,dotsize=6pt,linecolor=black](0,-2.5)
\rput[bl](-0.1,-3){\black{$\mathbf{e_2}$}}
\psdots[dotstyle=*,dotsize=6pt,linecolor=black](-2.165,1.25)
\rput[bl](-2.6,1.4){\black{$\mathbf{e_4}$}}
\psdots[dotstyle=*,dotsize=6pt,linecolor=black](0,0)
\rput[bl](0.2,0.5){\black{$\mathbf{e_1}$}}
\end{pspicture*}
\caption{Oriented Fano Plane:\\[8pt]
$P_{0}\,\,:   (a,b)(c,d)=(ca-b^*d,da^*+bc)$     $\downarrow$      $\circlearrowleft$
     $\rightarrow$\\
$P_{1}\,\,:   (a,b)(c,d)=(ca-db^*,a^*d+cb)$     $\downarrow $      $\circlearrowleft$
     $\leftarrow$\\
$P_{2}\,\,:   (a,b)(c,d)=(ac-b^*d,da^*+bc)$     $\downarrow $      $\circlearrowright$
     $\rightarrow$\\
$P_{3}\,\,:   (a,b)(c,d)=(ac-db^*,a^*d+cb)$     $\downarrow $      $\circlearrowright$
     $\leftarrow$\\
$\transp{P}_{0}:  (a,b)(c,d)=(ca-bd^*,ad+c^*b)$     $\uparrow$      $\circlearrowright$
     $\leftarrow$\\
$\transp{P}_{1}:  (a,b)(c,d)=(ca-d^*b,da+bc^*)$     $\uparrow$      $\circlearrowright$
     $\rightarrow$\\
$\transp{P}_{2}:  (a,b)(c,d)=(ac-bd^*,ad+c^*b)$     $\uparrow $      $\circlearrowleft$
     $\leftarrow$\\
$\transp{P}_{3}:  (a,b)(c,d)=(ac-d^*b,da+bc^*)$     $\uparrow $      $\circlearrowleft$
     $\rightarrow$
}
\label{fig: fano}
\end{figure}

An inspection of the eight products shows that for four of the products, the directional sense of the inner circle is the same as the directional sense of the three sides of the triangle. Either both are clockwise or both are counter-clockwise. For the other four, the senses are opposite, with one being clockwise and the other counter clockwise. Each one of the four products where the sense of the inner circle and the sides are the same produces an example of the octonions. The remaining four fail to satisfy the Moufang identities and have zero divisors at the third doubling. For those four, if $0\ne p\ne q\ne 0$ then $(e_{2p}+e_{2q+1})(e_{2p+1}+e_{2q})=0$. For example, for $p=1,q=2$ this gives $ (e_2+e_5)(e_3+e_4)=0 $. So strictly speaking, these four are not Cayley-Dickson doubling products. But they have been left in the discussion because of the interesting fractal patterns of their structure constants.

\section{The twists on W corresponding to the eight doubling products}

For each of the eight doubling products of finite sequences, its multiplication table is divided into distinct $2\times2$ blocks of the form

\[
\begin{bmatrix*}[l] e_{2r}e_{2s} & e_{2r}e_{2s+1}\\
                    e_{2r+1}e_{2s} & e_{2r+1}e_{2s+1}
\end{bmatrix*}
=
\begin{bmatrix*}[l] (e_{r},0)(e_{s},0)  & (e_{r},0)(0,e_{s})\\
                    (0,e_{r})(e_{s},0)  & (0,e_{r})(0,e_{s}) 
\end{bmatrix*}
\]

The results of the four products will vary according to which of the eight doubling products is used. Results are shown in Table \ref{table: blocks} on page \pageref{table: blocks}.

\begin{table}[ht]
\caption{Recursive Product of Units}\label{table: blocks}
\begin{tabular}{ll}

$\begin{array}{lll} P_{0}&e_{2s}&e_{2s+1}\\\toprule e_{2r}&(e_se_r,0)&(0,e_se_r^*)  \\\midrule e_{2r+1}  &(0,e_re_s)  &(-e_r^*e_s,0) \end{array}$
&
$\begin{array}{lll} \transp{P}_{0}&e_{2s}&e_{2s+1}\\\toprule e_{2r}&(e_se_r,0)&(0,e_re_s)  \\\midrule 
e_{2r+1}  &(0,e_s^*e_r)  &(-e_re_s^*,0) \end{array}$\\\bottomrule

$\begin{array}{lll} P_{1}&&\\\toprule e_{2r}&(e_se_r,0)&(0,e_r^*e_s)  \\\midrule e_{2r+1}  &(0,e_se_r)  &(-e_se_r^*,0) \end{array}$
&
$\begin{array}{lll} \transp{P}_{1}&&\\\toprule e_{2r}&(e_se_r,0)&(0,e_se_r)\\\midrule 
e_{2r+1}  &(0,e_re_s^*)&(-e_s^*e_r,0) \end{array}$\\\bottomrule

$\begin{array}{lll} P_{2}&&\\\toprule e_{2r}&(e_re_s,0)&(0,e_se_r^*)  \\\midrule e_{2r+1}  &(0,e_re_s)  &(-e_r^*e_s,0) \end{array}$
&
$\begin{array}{lll} \transp{P}_{2}&&\\\toprule e_{2r}&(e_re_s,0)&(0,e_re_s)  \\\midrule 
e_{2r+1}  &(0,e_s^*e_r)  &(-e_re_s^*,0) \end{array}$\\\bottomrule

$\begin{array}{lll} P_{3}&&\\\toprule e_{2r}&(e_re_s,0)&(0,e_r^*e_s)  \\\midrule e_{2r+1}  &(0,e_se_r) &(-e_se_r^*,0) \end{array}$
&
$\begin{array}{lll} \transp{P}_{3}&&\\\toprule e_{2r}&(e_re_s,0)&(0,e_se_r)  \\\midrule 
e_{2r+1}  &(0,e_re_s^*)  &(-e_s^*e_r,0) \end{array}$\\
\end{tabular}
\end{table}

The basis elements $e_0,e_1,e_2,\cdots$ are indexed by $W=\{0,1,2,\cdots\}$ which is a group under the bit-wise
`exclusive or' of their binary representations and for each of the 32 doubling products there is a function
$\omega$ from $W\times W$ to $\{-1,1\}$ such that for $p,q\in W$

\begin{equation}
  e_pe_q=\omega(p,q)e_{pq}
\end{equation}

where $pq$ is the group product of $p$ and $q$. The function $\omega$ is called a `twist' on the group $W$ and turns
the set of all finite sequences into a twisted group algebra.

Since $1=(1,0)=e_0$ is the identity, it follows that $\omega(p,0)=\omega(0,q)=1$ for all doubling products, and
since for $p>1,$ $e_pe_p^*=\|e_p\|^2=1$ and since $e_pe_p^*=-e_pe_p=-\omega(p,p)e_0=1$ it follows that for $p>1$

\begin{equation}
  \omega(p,p)=-1
\end{equation}

and that for all $p$

\begin{equation}
  e_p^*=\omega(p,p)e_p
\end{equation}

From equation \vref{eq:symmetric} it follows that

\begin{equation}\label{Eqn:neg}
  \omega(q,p)+\omega(p,q)=0\text{\ for\ }0\ne p\ne q\ne0
\end{equation}

\subsection{The product of $e_{2r}e_{2s}$}

Since $e_{2r}=\left(e_r,0\right)$ then $e_{2r}e_{2s}=\left(e_r,0\right)\left(e_s,0\right)$ this could be called the
case of $b=d=0.$ For each of the 32 distinct products, either $e_{2r}e_{2s}=\left(e_se_r,0\right)$ or
$e_{2r}e_{2s}=\left(e_re_s,0\right)$. We shall consider the effect each of these alternatives upon the twist $\omega$.

 \subsubsection{$e_{2r}e_{2s}=\left(e_se_r,0\right)$}
 
  Since $e_{2r}e_{2s}=\omega(2r,2s)e_{2rs}$ and since $e_{2r}e_{2s}=\omega(s,r)e_{2rs}$ we may conclude that
Whenever $e_{2r}e_{2s}=\left(e_se_r,0\right)$,
  \begin{equation}
     \omega(2r,2s)=\omega(s,r)
  \end{equation}

 \subsubsection{$e_{2r}e_{2s}=\left(e_re_s,0\right)$}
 
   Since $e_{2r}e_{2s}=\omega(2r,2s)e_{2rs}$ and since $e_{2r}e_{2s}=\omega(r,s)e_{2rs}$ we may conclude that
Whenever $e_{2r}e_{2s}=\left(e_re_s,0\right)$,
  \begin{equation}
     \omega(2r,2s)=\omega(r,s)
  \end{equation}

\subsection{The products $e_{2r}e_{2s+1}$ and $e_{2r+1}e_{2s}$}

Since $e_{2r}e_{2s+1}=\left(e_r,0\right)\left(0,e_s\right)$ this could be called the case of $b=c=0$. And since
$e_{2r+1}e_{2s}=\left(0,e_r\right)\left(e_s,0\right)$ that could be called the case of $a=d=0.$ An inspection of the 32
alternate products shows that there are only four distinct possibilities for the products of  $e_{2r}e_{2s+1}$ and
$e_{2r+1}e_{2s}$.

 \subsubsection{ $e_{2r}e_{2s+1 }=\left(0 ,e_se_r  \right)$ and $e_{2r+1 }e_{2s }=\left(0,e_re_s^* 
\right)$}
 
 These two conditions imply that whenever $e_{2r}e_{2s+1 }=\left(0 ,e_se_r  \right)$,
 
 \begin{equation}
    \omega(2r,2s+1)=\omega(s,r)
 \end{equation}
and that whenever $e_{2r+1 }e_{2s }=\left(e_re_s^* ,0  \right)$,
\begin{equation}
  \omega(2r+1,2s)=\omega(s,s)\omega(r,s)=\begin{cases}-\omega(r,s)\text{\ if\ }s>0\\
                                                      1\text{\ otherwise}
                                         \end{cases}
\end{equation}

 \subsubsection{$e_{2r}e_{2s+1 }=\left(0 ,e_re_s  \right)$ and $e_{2r+1 }e_{2s }=\left(0,e_s^*e_r 
\right)$}
 
 These two conditions imply that whenever $e_{2r}e_{2s+1 }=\left(0,e_re_s  \right)$,
 
 \begin{equation}
    \omega(2r,2s+1)=\omega(r,s)
 \end{equation}
and that whenever $e_{2r+1 }e_{2s }=\left(e_s^*e_r ,0  \right)$,
\begin{equation}
  \omega(2r+1,2s)=\omega(s,s)\omega(s,r)=\begin{cases}-\omega(s,r)\text{\ if\ }s>0\\
                                                      1\text{\ otherwise}
                                         \end{cases}
\end{equation}

 \subsubsection{ $e_{2r}e_{2s+1 }=\left(0 ,e_se_r^*  \right)$ and $e_{2r+1 }e_{2s }=\left(0,e_re_s 
\right)$}
 
These two conditions imply that whenever $e_{2r}e_{2s+1 }=\left(0 ,e_se_r^*  \right)$
 
 \begin{equation}
    \omega(2r,2s+1)=\omega(r,r)\omega(s,r)=\begin{cases}
                                             -\omega(s,r)\text{\ if\ }r>0\\
                                             1\text{\ otherwise}
                                           \end{cases}
 \end{equation}
 
and that whenever $e_{2r+1 }e_{2s }=\left(e_re_s ,0  \right)$,

\begin{equation}
  \omega(2r+1,2s)=\omega(r,s)
\end{equation}

 \subsubsection{$e_{2r}e_{2s+1 }=\left(0 ,e_r^*e_s  \right)$ and $e_{2r+1 }e_{2s }=\left(0,e_se_r 
\right)$}
 
 These two conditions imply that whenever $e_{2r}e_{2s+1 }=\left(0 ,e_r^*e_s  \right)$
 
 \begin{equation}
    \omega(2r,2s+1)=\omega(s,s)\omega(r,s)=\begin{cases}
                                             -\omega(r,s)\text{\ if\ }s>0\\
                                             1\text{\ otherwise}
                                           \end{cases}
 \end{equation}
 
and that whenever  $e_{2r+1 }e_{2s }=\left(e_se_r ,0  \right)$,

\begin{equation}
  \omega(2r+1,2s)=\omega(r,s)
\end{equation}

\subsection{The product $e_{2r+1}e_{2s+1}$}

Since $e_{2r+1}=\left(0,e_r\right)$ the product \\$e_{2r+1}e_{2s+1}$ could be called the case of $a=c=0.$

 \subsubsection{$e_{2r+1 }e_{2s+1 }=-\left(e_s^*e_r,0  \right)$}
 
 This implies that whenever $e_{2r+1 }i_{2s+1 }=-\left(e_s^*e_r,0  \right)$
 
 \begin{equation}
   \omega(2r+1,2s+1)=-\omega(s,s)\omega(s,r)=\begin{cases}
                                               \omega(s,r)\text{\ if\ }s>0\\
                                             -1\text{\ otherwise}  
                                             \end{cases}
 \end{equation}

 \subsubsection{$e_{2r+1 }e_{2s+1 }=-\left(e_re_s^*,0  \right)$}
 
 This implies that whenever $e_{2r+1 }e_{2s+1 }=-\left(e_re_s^*,0)  \right)$
 
  \begin{equation}
   \omega(2r+1,2s+1)=-\omega(s,s)\omega(r,s)=\begin{cases}
                                               \omega(r,s)\text{\ if\ }s>0\\
                                             -1\text{\ otherwise}  
                                             \end{cases}
 \end{equation}

 \subsubsection{$e_{2r+1 }e_{2s+1 }=-\left(e_se_r^*,0  \right)$}
 
 This implies that whenever $e_{2r+1 }e_{2s+1 }=-\left(e_se_r^*,0  \right)$
 
 \begin{equation}
   \omega(2r+1,2s+1)=-\omega(r,r)\omega(s,r)=\begin{cases}
                                               \omega(s,r)\text{\ if\ }r>0\\
                                             -1\text{\ otherwise}  
                                             \end{cases}
 \end{equation}

 \subsubsection{$e_{2r+1 }e_{2s+1 }=-\left(e_r^*e_s,0  \right)$}
 
This implies that whenever $e_{2r+1 }e_{2s+1 }=-\left(e_r^*e_s,0  \right)$

 \begin{equation}
   \omega(2r+1,2s+1)=-\omega(r,r)\omega(r,s)=\begin{cases}
                                               \omega(r,s)\text{\ if\ }r>0\\
                                             -1\text{\ otherwise}  
                                             \end{cases}
 \end{equation}
 
\section{The Eight Variations}
The preceding results are summarized in the following tables:

\begin{enumerate}
     \item[] \fbox{$P_{0}  :(a,b)(c,d)=(ca-b^*d,da^*+bc)$}\\  
\fbox{$\begin{array}{lll}\omega_{0}&2s&2s+1\\\toprule 2r&\omega(s,r)&\begin{cases}-\omega(s,r)\text{\ if\
}r>0\\1\text{\ otherwise}\end{cases}  \\\midrule 2r+1  
&\omega(r,s)  
&\begin{cases}\omega(r,s)\text{\ if\ }r>0\\-1\text{\ otherwise}\end{cases}\end{array}$}
    
    \item[] \fbox{$\transp{P}_{0}   :(a,b)(c,d)=(ca-bd^*,ad+c^*b)$}\\  
\fbox{$\begin{array}{lll}\omega^*_{0}&2s&2s+1\\\toprule 2r&\omega(s,r)&\omega(r,s)  \\\midrule 2r+1  
&\begin{cases}-\omega(s,r)\text{\ if\ }s>0\\ 1\text{\ otherwise}\end{cases}  
&\begin{cases}\omega(r,s)\text{\ if\ }s>0\\-1\text{\ otherwise}\end{cases}\end{array}$}

     \item[] \fbox{$P_{1}  :(a,b)(c,d)=(ca-db^*,a^*d+cb)$}\\  
\fbox{$\begin{array}{lll}\omega_{1}&2s&2s+1\\\toprule 2r&\omega(s,r)&\begin{cases}-\omega(r,s)\text{\ if\
}r>0\\1\text{\ otherwise}\end{cases}  \\\midrule 2r+1  
&\omega(s,r)  
&\begin{cases}\omega(s,r)\text{\ if\ }r>0\\-1\text{\ otherwise}\end{cases}\end{array}$}

     \item[] \fbox{$\transp{P}_{1}   :(a,b)(c,d)=(ca-d^*b,da+bc^*)$}\\  
\fbox{$\begin{array}{lll}\omega^*_{1}&2s&2s+1\\\toprule 2r&\omega(s,r)&\omega(s,r)\\\midrule 2r+1  
&\begin{cases}-\omega(r,s)\text{\ if\ }s>0\\ 1\text{\ otherwise}\end{cases}
&\begin{cases}\omega(s,r)\text{\ if\ }s>0\\-1\text{\ otherwise}\end{cases}\end{array}$}

     \item[] \fbox{$P_{2}  :(a,b)(c,d)=(ac-b^*d,da^*+bc)$}\\  
\fbox{$\begin{array}{lll}\omega_{2}&2s&2s+1\\\toprule 2r&\omega(r,s)&\begin{cases}-\omega(s,r)\text{\ if\
}r>0\\1\text{\ otherwise}\end{cases}  \\\midrule 2r+1  
&\omega(r,s)  
&\begin{cases}\omega(r,s)\text{\ if\ }r>0\\-1\text{\ otherwise}\end{cases}\end{array}$}

     \item[] \fbox{$\transp{P}_{2}  :(a,b)(c,d)=(ac-bd^*,ad+c^*b)$}\\  
\fbox{$\begin{array}{lll}\omega^*_{2}&2s&2s+1\\\toprule 2r&\omega(r,s)&\omega(r,s)  \\\midrule 2r+1  
&\begin{cases}-\omega(s,r)\text{\ if\ }s>0\\ 1\text{\ otherwise}\end{cases}  
&\begin{cases}\omega(r,s)\text{\ if\ }s>0\\-1\text{\ otherwise}\end{cases}\end{array}$}

     \item[] \fbox{$P_{3}  :(a,b)(c,d)=(ac-db^*,a^*d+cb)$}\\  
\fbox{$\begin{array}{lll}\omega_{3}&2s&2s+1\\\toprule 2r&\omega(r,s)&\begin{cases}-\omega(r,s)\text{\ if\
}r>0\\1\text{\ otherwise}\end{cases}  \\\midrule 2r+1  
&\omega(s,r) 
&\begin{cases}\omega(s,r)\text{\ if\ }r>0\\-1\text{\ otherwise}\end{cases}\end{array}$}

     \item[] \fbox{$\transp{P}_{3}  :(a,b)(c,d)=(ac-d^*b,da+bc^*)$}\\  
\fbox{$\begin{array}{lll}\omega^*_{3}&2s&2s+1\\\toprule 2r&\omega(r,s)&\omega(s,r)\\\midrule 2r+1  
&\begin{cases}-\omega(r,s)\text{\ if\ }s>0\\ 1\text{\ otherwise}\end{cases}  
&\begin{cases}\omega(s,r)\text{\ if\ }s>0\\-1\text{\ otherwise}\end{cases}\end{array}$}
\end{enumerate}

\section{The Twist Blocks}

The twist tables for each of the eight $\omega_k$ twist functions can each be subdivided into $2\times2$ blocks with entries $+1$ or $-1$. The upper left corner $C$ of each twist table for $\cda{N}$, $N>0$, corresponds to $r=s=0$. For $N>1$ the left side blocks $L$ correspond to $r>s=0$. The top side blocks $T$ correspond to $0=r<s.$ The diagonal blocks $-D$ correspond to $0<r=s.$ For $N>2$, the interior blocks $N$ correspond to $0\ne r \ne s\ne 0.$ Note that $D$ is shown in its positive form in Table \ref{Table: twistblocks} on page \pageref{Table: twistblocks} but appears as $-D$ in the actual twist tables. Also notice that $C$ is the same for all eight doubling products. Doubling products $P_0$ through $P_3$ differ from each other only in their values for $N$ and products $\transp{P}_0$ through $\transp{P}_3$ also differ from each other only in their values of $N.$ The twist tables for $\transp{P}_k$ is the transpose of the twist table for $P_k.$

Since every pair of basis vectors $e_p$, $e_q$ has a product $e_pe_q=\omega_k(p,q)e_{pq}$ where $pq=p\oplus q$, the bit-wise `exclusive or' (XOR) of the binary representations of $p$ and $q$ one only needs a table of values of $\omega_k(p,q)$ or $\omega_k^*(p,q)$ to recover the multiplication table for any particular $0\le k<4$ and any $\cda{N}$.

\small
\begin{table}[ht]
  \caption{Twist Blocks}
  \label{Table: twistblocks}
\begin{tabular}{cccccc}
 Twist  & $\mathbf{C}$  & $\mathbf{L}$ & $\mathbf{T}$   & $\mathbf{D}$   & $\mathbf{N}$ \\\toprule
        & $r=s=0$ & $r>s=0$   & $s>r=0$    & $r=s>0$  & $ 0\ne r\ne s\ne 0$ \\\midrule
 $\omega_{0}$ & $\left[ \begin{array}{rr}  1 &  1\\  1 & -1 \end{array} \right]$
              & $\left[ \begin{array}{rr}  1 & -1\\  1 &  1 \end{array} \right]$
              & $\left[ \begin{array}{rr}  1 &  1\\  1 & -1 \end{array} \right]$
              & $\left[ \begin{array}{rr}  1 & -1\\  1 &  1 \end{array} \right]$
              & $\left[ \begin{array}{rr} -1 &  \phantom{-}1\\  1 &  1 \end{array} \right]$\\\midrule

 $\omega_{1}$ & $\left[ \begin{array}{rr}  1 &  1\\  1 & -1 \end{array} \right]$
              & $\left[ \begin{array}{rr}  1 & -1\\  1 &  1 \end{array} \right]$
              & $\left[ \begin{array}{rr}  1 &  1\\  1 & -1 \end{array} \right]$
              & $\left[ \begin{array}{rr}  1 & -1\\  1 &  1 \end{array} \right]$
              & $\left[ \begin{array}{rr} -1 & -1\\ -1 & -1 \end{array} \right]$\\\midrule

 $\omega_{2}$ & $\left[ \begin{array}{rr}  1 &  1\\  1 & -1 \end{array} \right]$
              & $\left[ \begin{array}{rr}  1 & -1\\  1 &  1 \end{array} \right]$
              & $\left[ \begin{array}{rr}  1 &  1\\  1 & -1 \end{array} \right]$
              & $\left[ \begin{array}{rr}  1 & -1\\  1 &  1 \end{array} \right]$
              & $\left[ \begin{array}{rr}  \phantom{-}1 &  \phantom{-}1\\  1 &  1 \end{array} \right]$\\\midrule

 $\omega_{3}$ & $\left[ \begin{array}{rr}  1 &  1\\  1 & -1 \end{array} \right]$
              & $\left[ \begin{array}{rr}  1 & -1\\  1 &  1 \end{array} \right]$
              & $\left[ \begin{array}{rr}  1 &  1\\  1 & -1 \end{array} \right]$
              & $\left[ \begin{array}{rr}  1 & -1\\  1 &  1 \end{array} \right]$
              & $\left[ \begin{array}{rr}  1 & -1\\ -1 & -1 \end{array} \right]$\\\midrule

$\omega_{0}^*$ & $\left[ \begin{array}{rr}  1 &  1\\  1 & -1 \end{array} \right]$
              & $\left[ \begin{array}{rr}  1 &  1\\  1 & -1 \end{array} \right]$
              & $\left[ \begin{array}{rr}  1 &  1\\ -1 &  1 \end{array} \right]$
              & $\left[ \begin{array}{rr}  1 &  1\\ -1 &  1 \end{array} \right]$
              & $\left[ \begin{array}{rr} -1 &  \phantom{-}1\\  1 &  1 \end{array} \right]$\\\midrule

$\omega_{1}^*$ & $\left[ \begin{array}{rr}  1 &  1\\  1 & -1 \end{array} \right]$
              & $\left[ \begin{array}{rr}  1 &  1\\  1 & -1 \end{array} \right]$
              & $\left[ \begin{array}{rr}  1 &  1\\ -1 &  1 \end{array} \right]$
              & $\left[ \begin{array}{rr}  1 &  1\\ -1 &  1 \end{array} \right]$
              & $\left[ \begin{array}{rr} -1 & -1\\ -1 & -1 \end{array} \right]$\\\midrule
$\omega_{2}^*$ & $\left[ \begin{array}{rr}  1 &  1\\  1 & -1 \end{array} \right]$
              & $\left[ \begin{array}{rr}  1 &  1\\  1 & -1 \end{array} \right]$
              & $\left[ \begin{array}{rr}  1 &  1\\ -1 &  1 \end{array} \right]$
              & $\left[ \begin{array}{rr}  1 &  1\\ -1 &  1 \end{array} \right]$
              & $\left[ \begin{array}{rr}  \phantom{-}1 &  \phantom{-}1\\  1 &  1 \end{array} \right]$\\\midrule

$\omega_{3}^*$ & $\left[ \begin{array}{rr}  1 &  1\\  1 & -1 \end{array} \right]$
              & $\left[ \begin{array}{rr}  1 &  1\\  1 & -1 \end{array} \right]$
              & $\left[ \begin{array}{rr}  1 &  1\\ -1 &  1 \end{array} \right]$
              & $\left[ \begin{array}{rr}  1 &  1\\ -1 &  1 \end{array} \right]$
              & $\left[ \begin{array}{rr}  1 & -1\\ -1 & -1 \end{array} \right]$\\\midrule
\end{tabular}\\[12pt]
\end{table}
\normalsize

\section{Geometry of the $\omega$ twist tables}

The $\omega$ tables for the eight doubling products are highly periodic when using the shuffle basis. However, the two most common doubling products $\omega_3$ and $\omega_3^*$ have a somewhat bewildering pattern to the eye. This is not the case for $\omega_1$ and $\omega_2$ or their transposes. The pattern for $\omega_1$ `alternates' whereas as $N\to\infty$, the $\omega_2$ table for $\cda{N+1}$ has the same appearance as the table for $\cda{N}$. Note that for any $\omega$ table the table for $\cda{N}$ always forms the upper left quadrant of the table for $\cda{N+1}$.

\begin{figure}[ht]
\begin{multicols}{2}
\includegraphics[scale=0.3525]{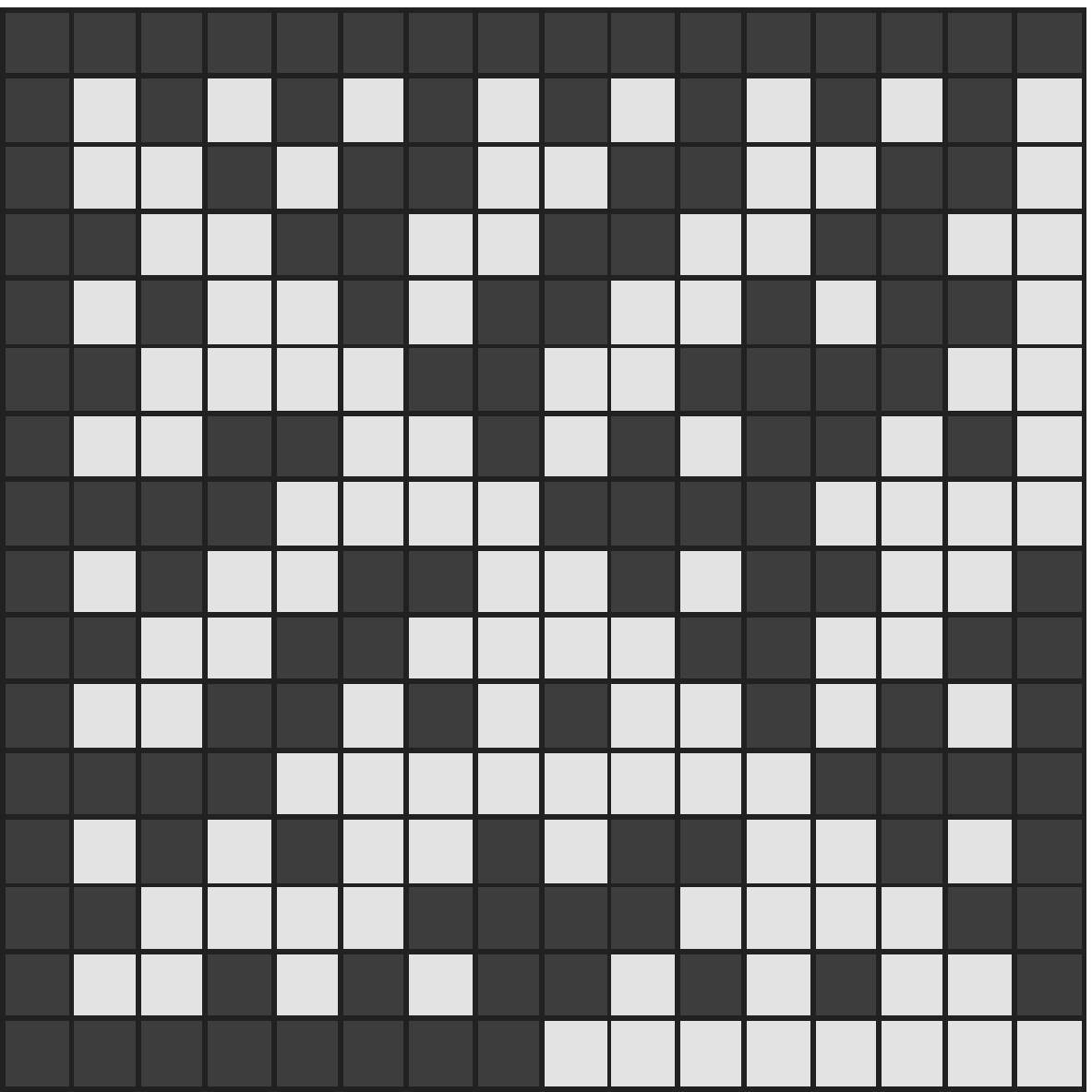}
\vfill
\columnbreak
\includegraphics[scale=0.354]{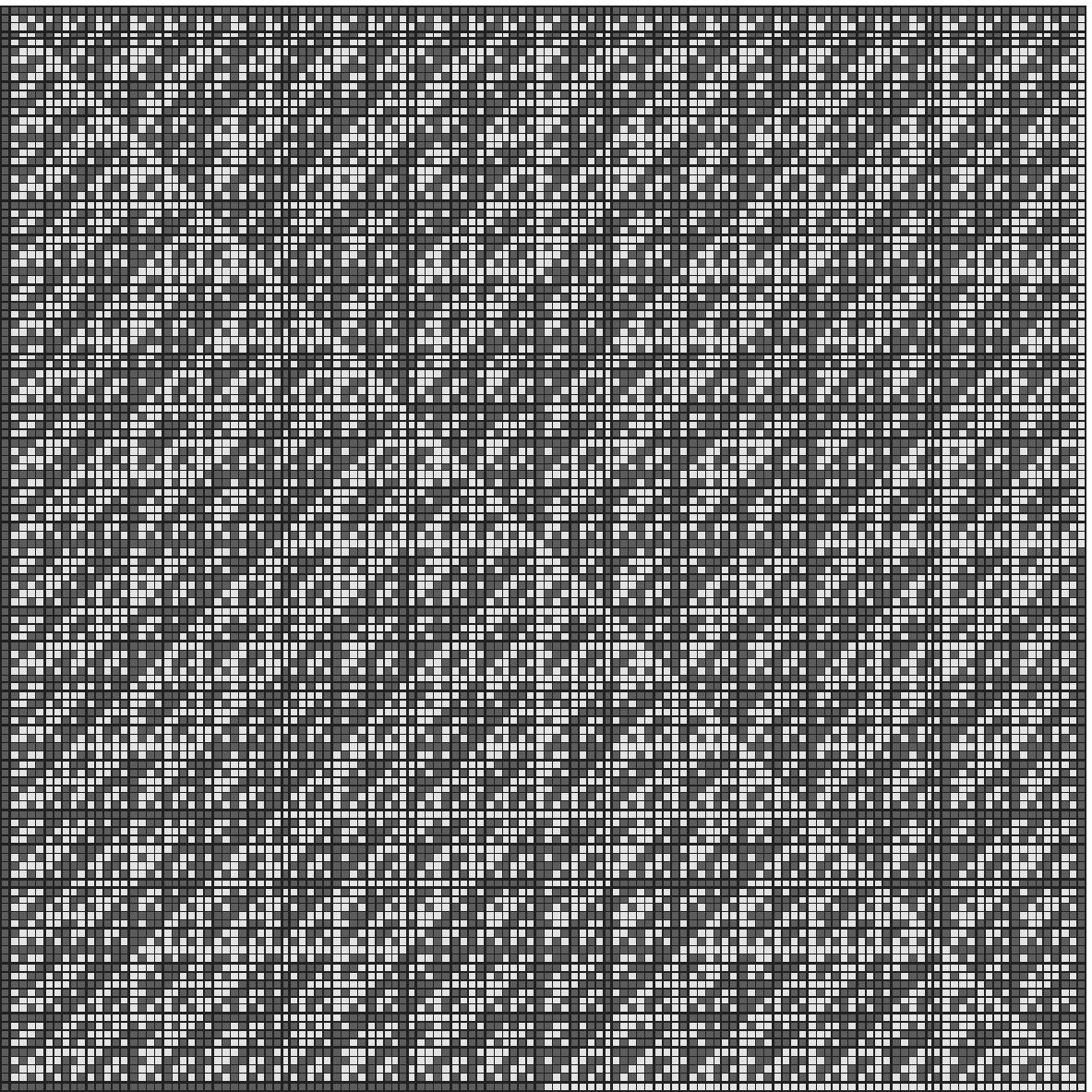}
\end{multicols}
\caption{ $\omega_0$ twist map for $\cda{4}$ and $\cda{7}$}\label{fig: omega0}
\end{figure}

\begin{figure}[ht]
\begin{multicols}{2}
\includegraphics[scale=0.3425]{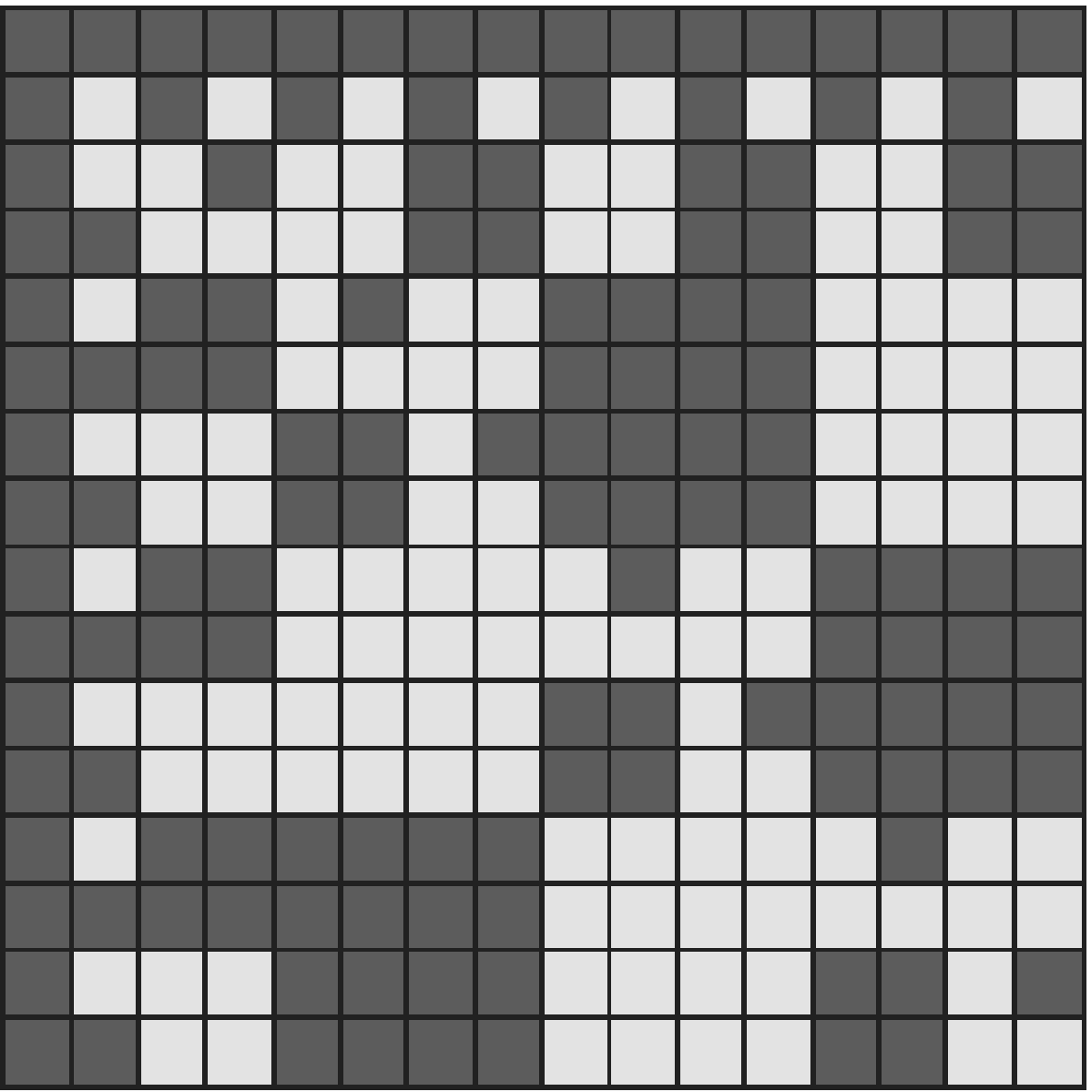}
\vfill
\columnbreak
\includegraphics[scale=0.354375]{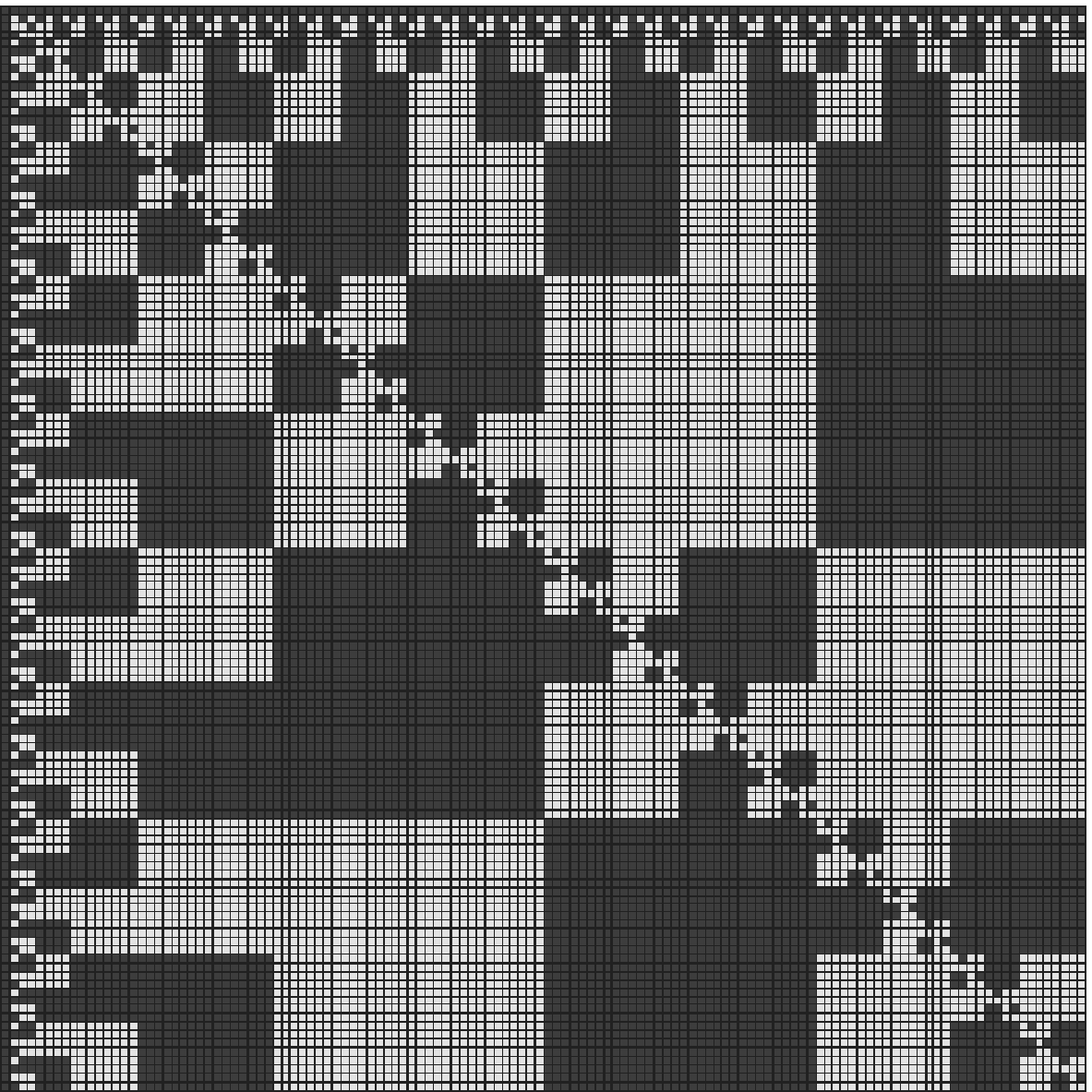}
\end{multicols}
\caption{ $\omega_1$ twist map for $\cda{4}$ and $\cda{7}$}\label{fig: omega1}
\end{figure}

\begin{figure}[ht]
\begin{multicols}{2}
\includegraphics[scale=0.3525]{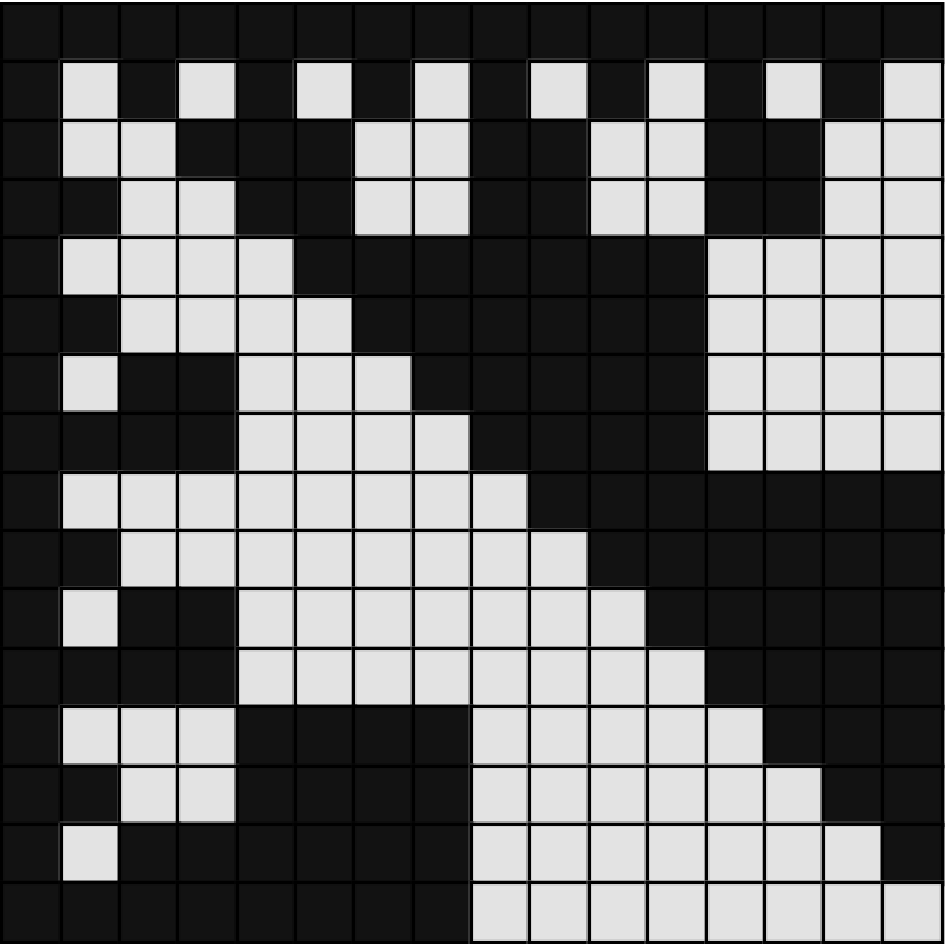}
\vfill
\columnbreak
\includegraphics[scale=0.28125]{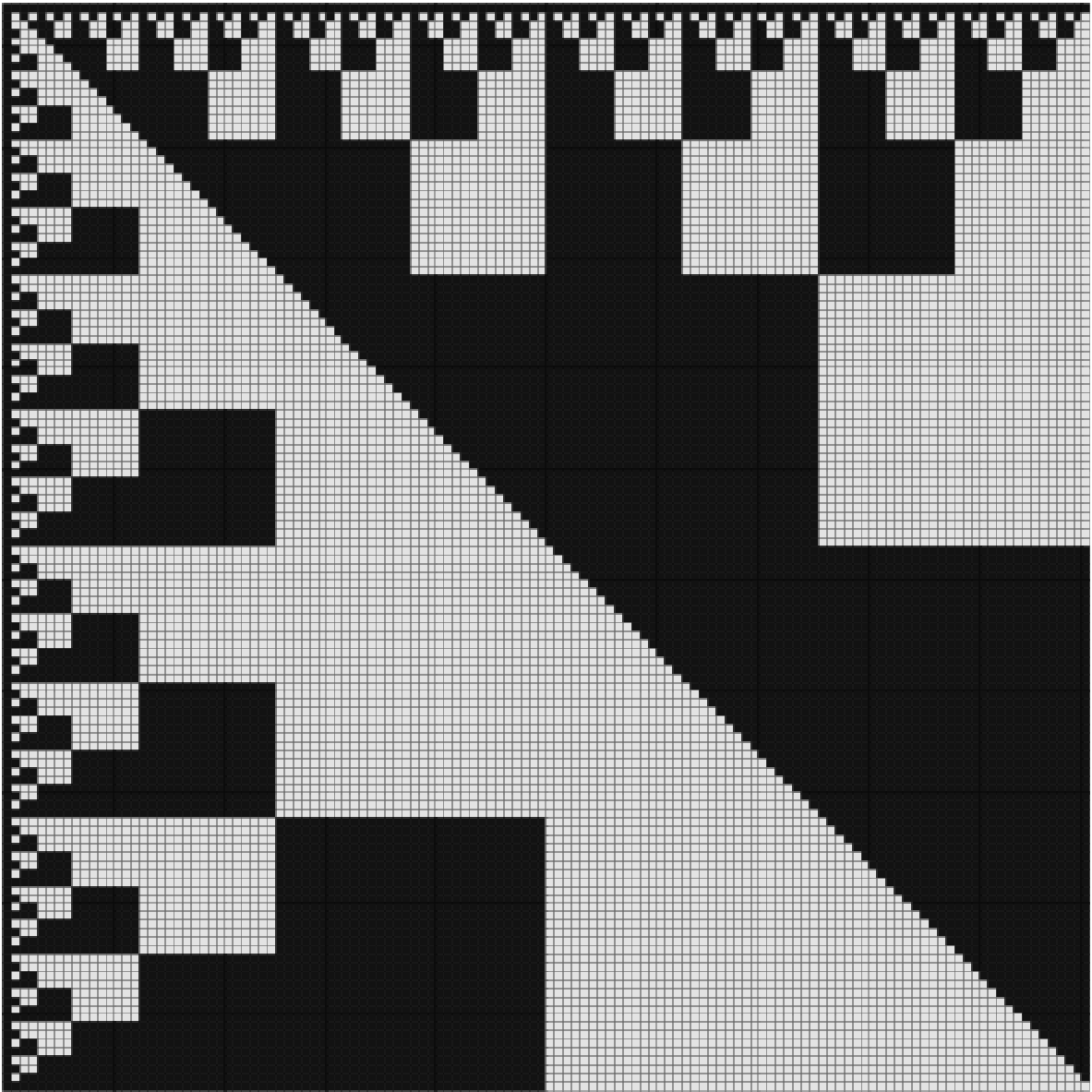}
\end{multicols}
\caption{ $\omega_2$ twist map for $\cda{4}$ and $\cda{7}$}\label{fig: omega2}
\end{figure}

\begin{figure}[ht]
\begin{multicols}{2}
\includegraphics[scale=0.3525]{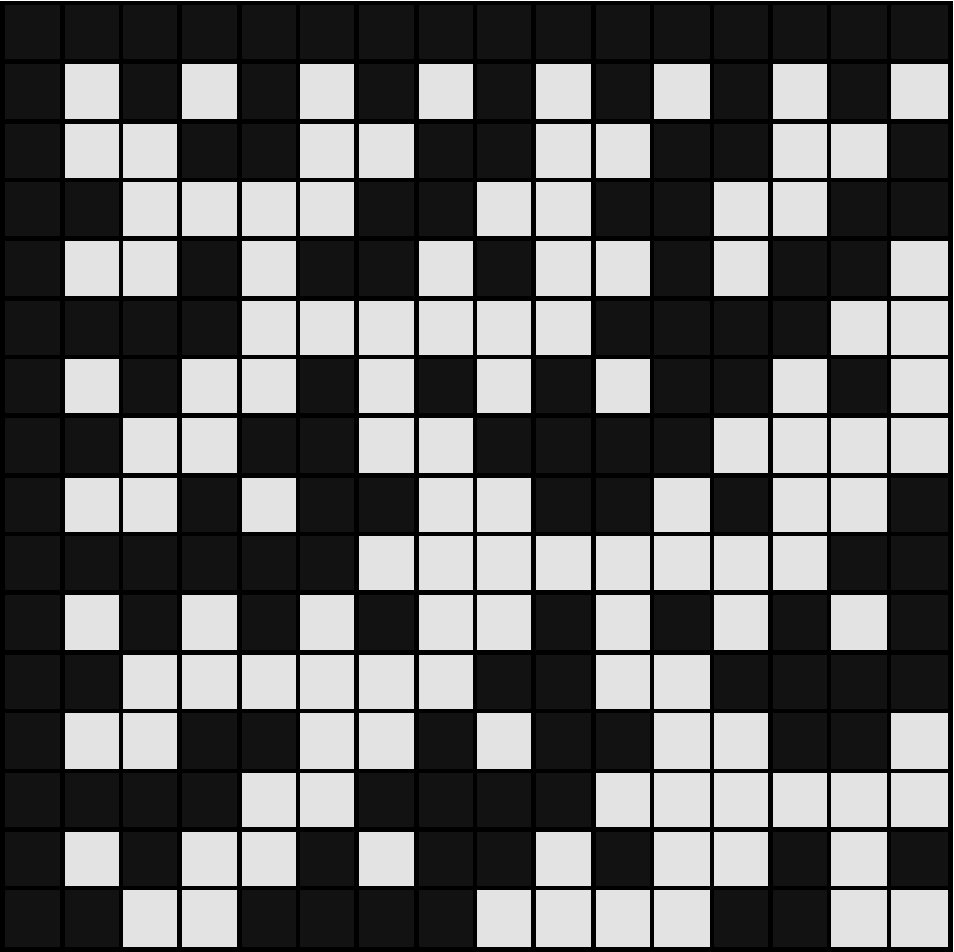}
\vfill
\columnbreak
\includegraphics[scale=0.354375]{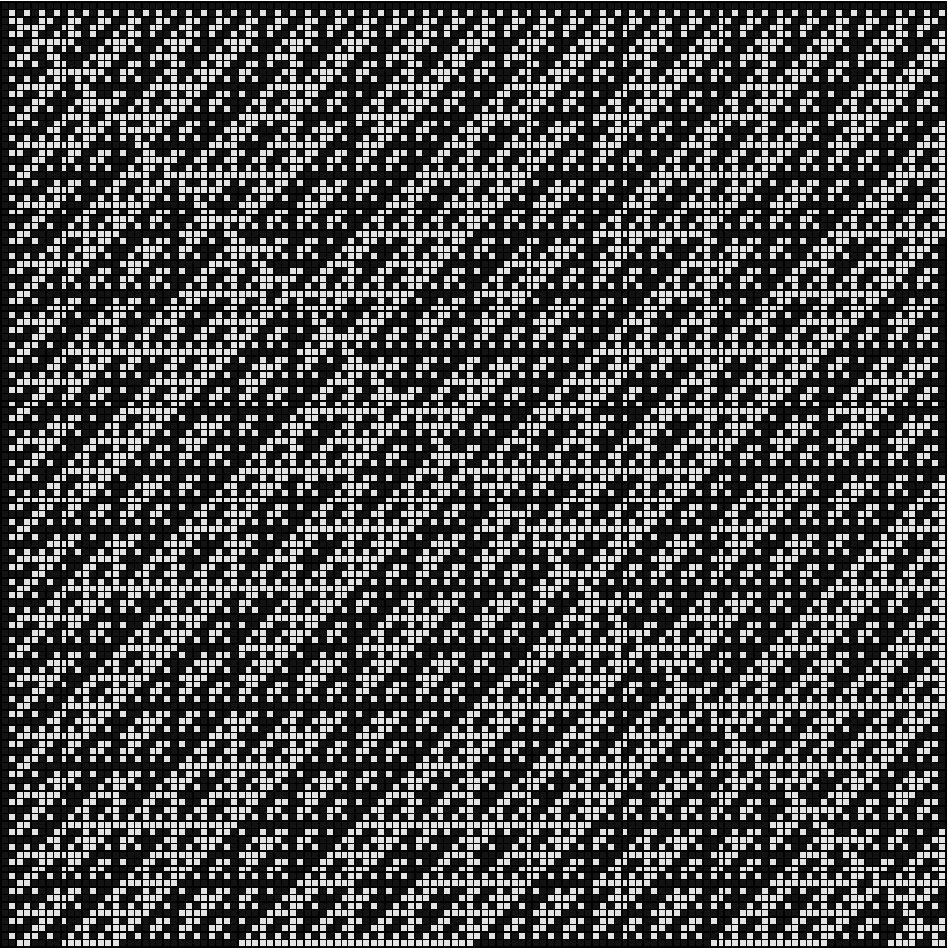}
\end{multicols}
\caption{ $\omega_3$ twist map for $\cda{4}$ and $\cda{7}$}\label{fig: omega3}
\end{figure}

In Figures \ref{fig: omega0}--\ref{fig: omega3} we see visual representations of $\omega_0$ through $\omega_3$ for $\cda{4}$ and $\cda{7}$ respectively, with dark gray representing $\omega(p,q)=+1$ and light gray representing $\omega(p,q)=-1$. The upper left square in all images has coordinates $(p,q)=(0,0)$ and the lower right corner has coordinates $(p,q)=\left(2^N-1,2^N-1\right)$.

\section{The twist tree for $\omega_2$}

In \cite{B2009} it was demonstrated how to use a tree diagram to determine $\omega_3(p,q)$ for any two non-negative integers $p$ and $q$. Figure \ref{fig:cydtree} is the \emph{twist tree} for $\omega_2$.

To use the tree to find $\omega_2(p,q)$ for non-negative integers $p$ and $q$ they must first be represented in their binary form. To illustrate, let $p=93=1011101_B$ and let $q=37=0100101_B$. Next `shuffle' the two bit strings together to obtain $10,01,10,10,11,00,11$. Beginning with the top of the tree $C$ which represents the top left quadrant of the $\omega_2$ twist table a binary 1 is an instruction to move to the right branch of the tree and a binary 0 is an instruction to move to the left branch of the tree. So the first instruction 10 moves from $C$ to $L$. The next instruction 01 moves from $L$ to $-1$. Therefore $\omega_2(93,37)=-1$. There is no need to finish traversing the sequence $10,01,10,10,11,00,11$ further because of the stability of $\omega_2$. Thus $e_{93}e_{37}=-e_{120}$ and $(37,93,120)$ is a structure constant for the product $P_2$.

Note that any sequence of instructions which terminates in $C,T$ or $L$ gives a result of $\omega_2(p,q)=1$ and any which terminates in $-D$ gives $\omega_2(p,q)=-1$.

This approach to finding the value of the twist function developed in \cite{B2009} was found useful by Flaut and Shpakivskyi in their study \cite{FS2015} of holomorphic functions in generalized Cayley-Dickson algebras.

\begin{table}[h]
~\\
\Tree [.$\mathbf{C}$ [. $\mathbf{C}$ $\mathbf{T}$ ] [. $\mathbf{L}$ $\mathbf{-D}$ ] ]
\Tree [.$\mathbf{L}$ [. $\mathbf{L}$ $-1$ ] [. $\mathbf{L}$ $+1$ ] ]\\[12pt]
\Tree [.$\mathbf{T}$ [. $\mathbf{T}$ $\mathbf{T}$ ] [. $+1$ $-1$ ] ]
\Tree [.$\mathbf{-D}$ [. $\mathbf{-D}$ $+1$ ] [ $-1$ $\mathbf{-D}$ ] ]
\vspace{3mm}
\caption{$P_2$ Twist Tree $(a,b)(c,d)=(ac-b^*d,da^*+bc)$}
\label{fig:cydtree}
\end{table}

The two dimensional graphs of $\omega_2$ in Figure \ref{fig: omega2} suggests that

\begin{equation}
 \omega_2(p,q) = \begin{cases}\phantom{-}1 \text{ if } 1<\frac{q}{p}<\frac{3}{2} \text{ and}\\
                                        -1 \text{ if } \frac{2}{3}<\frac{q}{p}\le 1
                 \end{cases}
\end{equation}

The trees in Table \ref{fig:cydtree} suggest that if $m>n\ge0$ and $r,s<2^n$ then

\begin{align}
 \omega_2\left(2^n+r,2^m+s\right)    &=1\\
 \omega_2\left(2^m+2^n+r,2^n+s\right)&=1\\
 \omega_2\left(2^m+r,2^m+2^n+s\right)&=1
\end{align}

\section{Conclusion}

Four variations of the doubling-products which produce the Cayley-Dickson algebras have been identified as well as four `pseudo' Cayley-Dickson algebras with fractal patterns in their structure constants. Also,  an alternate conception has been illustrated of the ordered pair of two sequences as consisting of the `shuffling' of the two sequences. We have also presented another instance of the use of a `twist tree' for determining the twist function of a Cayley-Dickson algebra. In regard to octonion algebras in particular we have suggested an alternate way to label the Fano Plane in such a way that the structure constants have a consistent orientation in the plane. It is hoped that the alternate perspectives offered in this paper will prove useful to other researchers.

\section*{Appendix}

The author wrote the following bc \cite{N2001} program for generating the $\omega$ twist maps in Figures \ref{fig: omega0} through \ref{fig: omega3}. The code generates a \LaTeX\ document containing a PSTricks figure. 

For $\cda{7}$ it was necessary to increase the LaTeX memory allotment in the user's (not the system's) texmf.cnf file to `\verb|main_memory = 7500000|'.

\small
\begin{verbatim}
#TeXTables.bc
define sgn(x,y){
 auto p,q,lp,lq;
#Note: lp means p is odd, !lp means p is even
 scale=0;
 p=x;
 q=y;
 if(p*q==0)             return 1;
 lp=p%2;
 lq=q%2;
  if(p==1 && lq)         return -1;
/*The following lines are for omega 0
  if(!lp && !lq) return  sgn(q/2,p/2);
  if(!lp &&  lq) return -sgn(q/2,p/2);
  return sgn(p/2,q/2);
End of omega 0 lines*/
/*The following lines are for omega 1
  if(!lp && lq) return  -sgn(p/2,q/2);
  return sgn(q/2,p/2);
End of omega 1 lines*/
#The following lines are for omega 2
  if(!lp && lq) return  -sgn(q/2,p/2);
  return sgn(p/2,q/2);
#End of omega 2 lines*/
/*The following lines are for omega 3
  if(!lp && !lq) return   sgn(p/2,q/2);
  if(!lp &&  lq) return  -sgn(p/2,q/2);
  return sgn(q/2,p/2);
End of omega 3 lines*/
}
k=7; #Table for Cayley-Dickson Algebra A_k (set desired value)
print "\\documentclass{minimal}\n";
print "\\usepackage{pstricks-ad
  if(p==1 && lq)         return -1;
/*The following lines are for omega 0
  if(!lp && !lq) return  sgn(q/2,p/2);
  if(!lp &&  lq) return -sgn(q/2,p/2);
  return sgn(p/2,q/2);
End of omega 0 lines*/
/*The following lines are for omega 1
  if(!lp && lq) return  -sgn(p/2,q/2);
  return sgn(q/2,p/2);
End of omega 1 lines*/
#The following lines are for omega 2
  if(!lp && lq) return  -sgn(q/2,p/2);
  return sgn(p/2,q/2);
#End of omega 2 lines*/
print "\\usepackage{graphicx}\n";
print "\\begin{document}\n";
print "\\scalebox{0.8}{\n";
print "\\begin{pspicture}(-13.5,-12.5)(7.5,12.5)\n";
scale=5;
 dot=24/2^k; # Sets the size of each square
 du=20/2^k;  # Sets the distance between squares
 v=12+du;   # Initializes y coordinate of square
for(p=0;p<2^k;p++){
 u=-13-du;  # Initializes x coordinate of square
 v-=du;     # Increments y coordinate of square
 for(q=0;q<2^k;q++){
    u+=du;  # Increments x coordinate of square
    if( sgn(p,q)==1  ) print " \\psdot[dotsize=",dot,
    "\\psxunit,dotstyle=square,fillcolor=blue]\n
    (",u,",",v,")\n" else print " \\psdot[dotsize=",dot,"
    \\psxunit,dotstyle=square,fillcolor=yellow]\n(",u,",",v,")
    \n"; /* remove last 4 carriage returns before executing */
 }
    print "\n";
}
print "\\end{pspicture}\n";
print "}\n";
print "\\end{document}";
quit
\end{verbatim}

\end{document}